\documentclass[final,leqno,onefignum,onetabnum]{siamltex1213}

\usepackage{bbm}
\usepackage{amsmath,amssymb}
\newtheorem{remark}[theorem]{Remark}
\newtheorem{assumption}[theorem]{Assumption}

\title{Observers for compressible Navier-Stokes equation}

\author{Amit Apte\footnotemark[1] \and Didier Auroux\footnotemark[2] \and Mythily Ramaswamy\footnotemark[3] \footnotemark[4]}

\begin{document}
\maketitle

\renewcommand{\thefootnote}{\fnsymbol{footnote}}

\footnotetext[1]{International Centre for Theoretical Sciences - TIFR, IISc Campus, Bangalore-560012, India. (\email{apte@icts.res.in})}
\footnotetext[2]{Universit\'e C\^ote d'Azur, CNRS, LJAD, France. (\email{auroux@unice.fr})}
\footnotetext[3]{T.I.F.R Centre for Applicable Mathematics, Post Bag No. 6503, GKVK Post Office, Bangalore-560065, India. (\email{mythily@math.tifrbng.res.in})}

\footnotetext[4]{The authors would like to thank the Indo-French
  Centre for Applied Mathematics (IFCAM) for financial support under
  the ``Observers and data assimilation'' project.}

\renewcommand{\thefootnote}{\arabic{footnote}}

\begin{abstract}
  We consider a multi-dimensional model of a compressible fluid in a
  bounded domain. We want to estimate the density and velocity of the
  fluid, based on the observations for only velocity. We build an
  observer exploiting the symmetries of the fluid dynamics laws. Our
  main result is that for the linearised system with full observations
  of the velocity field, we can find an observer which converges to
  the true state of the system at any desired convergence rate for
  finitely many but arbitrarily large number of Fourier modes.  Our
  one-dimensional numerical results corroborate the results for the
  linearised, fully observed system, and also show similar convergence
  for the full nonlinear system and also for the case when the
  velocity field is observed only over a subdomain.
\end{abstract}

\begin{keywords} Data assimilation; Observer; Navier-Stokes
  equation; \end{keywords}

\begin{AMS}  93C20; 93C95; 93B40 \end{AMS}

\pagestyle{myheadings}
\thispagestyle{plain}

\section{Introduction} \label{sec-intro}

Data assimilation is the problem of estimating the state of a
dynamical system described by an evolution equation, typically partial
differential equations (PDE), using observations, often noisy and
partial, of that system. This has been widely studied in the
geophysical context, e.g.  meteorology, oceanography, fluid flows,
etc., \cite{Ben01, kalnaybook, book:lewisvarahan06}. One of the
approaches to this state estimation problem, and the one that we study
in this paper, is the construction of appropriate observers
\cite[Chap.~7]{sontagbook}.

Observers are essentially a modification of the original evolution
equation for the system, to incorporate the observations in a feedback
term, with the aim that the solution of the observer converges to the
solution of the original system being observed. Observers for finite
dimensional systems have been well studied in the literature, see, for
example \cite{sontagbook, zabczykbook}. But in many applications such
as earth sciences or engineering, the systems are modeled using PDE
that are highly nonlinear and, in many instances, chaotic. In such
cases of infinite dimensional systems governed by PDE, there are only
a few examples available in the literature, mostly for linear systems
and very few for nonlinear systems \cite{lasieckaTriggiani,
  SmyshlyaevKrstic, XuDeguenon, deguenon2, GuoShao, GuoXu,
  johna:PODobs, vazquezMHD, mohamed14:observer, BensoussanEtAl,
  ramdaniTW10}.

Some of the commonly used observers are Kalman Filters or Luenberger
observers but the main drawback of these is that they often may break
intrinsic properties of the model, e.g. symmetries and/or physical
constraints such as balances in geophysical models, see for example
\cite{gottwald14:balance} and references therein.  For nonlinear
system possessing certain symmetries, it is natural to seek a
correction term which also preserves those symmetries. Such an
invariance may make the correction term nonlocal but it may have other
desirable properties, for example, independence from the change of
coordinates. Recently, there have been attempts to construct, for a
variety of systems, observers based on considerations of symmetry
\cite{bonnabelMR08, bonnabelMR09, auroux-bonnabel11}. The observers we
construct are motivated by these recent works, as we will see in
section~\ref{sec-ns-observer}. Another motivation for observer design
is the computational cost: a full Kalman filter is usually too
expensive for real applications, due to the size of the gain matrices,
while an observer might be much more affordable, without degrading the
identification process.

The main aim of this work is to develop an observer for a class of PDE
inspired by earth system applications, all of which use some
approximations of Navier-Stokes equations for fluid flow. In
particular, we propose an observer for the nonlinear PDE
[equations~\eqref{eqrhou}] describing the evolution of a compressible,
adiabatic fluid whose velocity field is observed either fully or
partially. This is justified by the fact that velocity observations
for fluids have become more available in the last years. We can cite
floats (e.g. Argo) for oceanography, balloons for meteorology, but
also all observations generated by optical flow or particle imaging
velocimetry from images (including satellite images). But as we will
see, the other case in which the density is observed and not the
velocity is very similar. Apart from a purely theoretical motivation,
the observations of sea surface density are now available \cite[and
references therein]{droghei2016combining} and thus the study of
observers with density observations is relevant.

Our main theoretical result, theorem~\ref{mainprop}, supported by
extensive numerical results of the one dimensional system
(section~\ref{sec-numeric}), is that in the case of complete
observations of the velocity, an appropriate choice of parameters
leads to convergence of the observer to the true solution at any
prescribed rate for arbitrarily large but finite number of Fourier
modes.

There has been some work on developing observers for Navier-Stokes
based systems. The papers \cite{vazquezMHD, johna:PODobs,
  parkguay:estimation} work with finite dimensional approximations of
the PDE involved whereas we work with the full PDE itself.  A
completely different approach based on developing observers using
appropriate continuous time limits of discrete time 3D-var or Kalman
filter is developed in a series of papers \cite{kls14, lss14, blsz13}
and also in \cite{bessaih:contDA, albanez:contDA, farhat:contDA,
  olsontiti03, hayden2011:discreteda, azouani:continuousDA}. All these
papers are in the context of incompressible flows. Here we deal with
the compressible model, a coupled system for density and velocity. The
linearized system consists of a hyperbolic and parabolic PDE.  This
coupling of mixed types poses some difficulties in tackling the system
theoretically, in particular the full nonlinear system, unlike in the
case of incompressible Navier-Stokes system. We refer to \cite{Erve}
and \cite{Chow}, where these difficulties have been overcome for the
nonlinear system, with highly involved techniques.  In the case of the
linearized system around constant steady states, we manage using
Fourier series. One of the main contributions of this paper is that we
also derive the decay rates for the convergence of the observers and
indeed find observers that can decay arbitrarily fast. Since
Navier-Stokes based PDE are commonly used in practical data
assimilation problems in earth sciences, our work has the potential to
be directly relevant to these applications, as we discuss in
section~\ref{sec-conclude}.

The outline of the paper is as follows. The multi-dimensional model
and the observer are both introduced in the next
section~\ref{sec-ns-observer}. In that section, we also state the
linearized version of this problem. We analyze the convergence of this
proposed observer for the linearized PDE in section~\ref{sec-theory}
and prove the main result theorem~\ref{mainprop}. We also briefly
discuss the difficulties that arise in the theoretical analysis of the
nonlinear equation or of the cases with either partial observations or
unknown forcing for the linearized equations. In
section~\ref{sec-numeric}, we present one-dimensional numerical
results to substantiate the linear theory. We also present the
numerical results showing the efficacy of the observer in estimating
the true solution for the partially observed linear system, the
observer in the case of unknown forcing term, and the fully nonlinear
case. The last section~\ref{sec-conclude} discusses some future
directions of research.

The authors would like to thank both the anonymous referees for useful
remarks which led to this improved revised version.

\section{Compressible Navier-Stokes equations and
  observers}\label{sec-ns-observer}

In this paper, we study a model of compressible fluid in a bounded
domain. The density $\rho(t,x)$ and velocity $u(t,x)$ of the fluid
form the state vector of this model and they obey the following
compressible Navier-Stokes system in $n$-dimensions:
\begin{eqnarray}
  \rho_t + \nabla \cdot (\rho u) = 0 \,, \quad 
  \rho \left[ u_t + (u \cdot \nabla) u \right] = -\nabla p(\rho) + \mu
  \Delta u + (\lambda + \mu) \nabla (\nabla \cdot u) \,, \label{eqrhou}
\end{eqnarray}
where $\lambda$ and $\mu$ are the Lam\'e parameters, satisfying the
standard assumptions of $\mu > 0$ and $\lambda + 2\mu/3 \ge 0$
\cite{valli1986navier, vaigant1995existence}. The pressure $p$ is
given by the adiabatic equation of state, $p(\rho) = \rho^\gamma$ with
$\gamma$, the adiabatic exponent, taken to be $1.4$ for the numerical
results discussed in section~\ref{sec-numeric}. We consider solutions
over a finite time interval $[0,T]$ for the space domain
$[0,1]^n \subset \mathbb{R}^n$ with periodic boundary conditions:
\begin{equation}
  \rho(t,x)=\rho(t,x+e_k)\,,\quad \forall\
  (t,x) \in [0,T] \times [0,1]^n \,, \quad \forall\ 1 \le k
  \le n \,, \label{eqbc}
\end{equation}
and similar periodic conditions for $u(t,x)$ and for all derivatives
of $\rho$ and $u$, where $e_k = (0, \dots, 1, \dots, 0)$ with the $1$
in the $k$-th place.

We assume that the initial conditions $\rho(0,x) = \rho_I(x)$ and
$u(0,x) = u_I(x)$ are unknown, but we have some information on the
solution $u(t,x)$ of equations~\eqref{eqrhou}. The goal is to build an
observer $(\hat\rho,\hat u)$ for this system, in such a way that both
the density and the velocity of the observer $(\hat\rho,\hat u)$
converges towards the solution $(\rho,u)$. We will assume that we have
observations of the velocity denoted by $u(t,x)$ for $t \in [0,T]$ and
$x \in \Omega \subset [0,1]^n$ where $\Omega$ is a subset of
$[0,1]^n$.  For most of the theoretical study we will consider
$\Omega = [0,1]^n$ while for the numerical study of the
one-dimensional problem, we will consider $\Omega = [0,L]$ with
$0 < L \le 1$.

\subsection{Observers} \label{subsec-observer}

We introduce the observer $(\hat\rho, \hat u)$, based on the
symmetries for the system, satisfying the following set of equations:
\begin{equation}
  \begin{array}{l}
    \hat\rho_t + \nabla \cdot (\hat\rho \hat u) = F_\rho(\hat\rho,\hat
    u,u)  \,, \\[0.1cm] \hat\rho \left[ \hat u_t + (\hat u \cdot
    \nabla) \hat u \right] = -\nabla p(\hat\rho) + \mu  \Delta \hat u
    + (\lambda + \mu) \nabla (\nabla \cdot \hat u) + F_u(\hat\rho,\hat
    u,u) \,.  \label{eqobsrhou}
 \end{array}
\end{equation}
As there are no observations of the density, we first assume that the
feedback terms do not depend on $\hat\rho$. Additionally, since these
terms should be equal to $0$ when $u$ and $\hat u$ coincide, it is
reasonable to consider the following classes of functions:
\begin{eqnarray}
  F_\rho(\hat\rho,\hat u,u) = \varphi_\rho\ast D_\rho(u-\hat
  u)\,, \quad F_u(\hat\rho,\hat u,u) = \varphi_u\ast
  D_u(u-\hat u),
\end{eqnarray}
where $D_\rho$ and $D_u$ are differential (or integral) operators
(e.g. $\nabla$, $\partial_t$, \dots), and $\varphi_\rho$ and
$\varphi_u$ are convolution kernels. We assume that these kernels are
time independent and hence only functions of $x$, and we assume that
they are isotropic. With this choice, the observer preserves the
symmetries of the system, since the correction terms are based on a
convolution product with an isotropic kernel, making it invariant by
rotation or translation. Moreover, observations may be noisy, so the
involvement of a convolution kernel leads to an observer that is
inherently robust to noise.

Note that for the case of observations over partial domain, i.e., when
$\Omega \ne [0,1]^n$, the feedback terms are present only in $\Omega$
so that we essentially write the feedback term as:
\begin{align}
  F_\rho(\hat\rho,\hat u,u) = \varphi_\rho\ast D_\rho \left[
     (u-\hat u) \mathbbm{1}_\Omega \right] \,, \quad F_u(\hat\rho,\hat
  u,u) = \varphi_u\ast D_u\left[ (u-\hat u) \mathbbm{1}_\Omega
    \right] \,. \label{feedbackrhou}
\end{align}

We would like to examine the convergence of the observer solution
$(\hat\rho, \hat u)$ to the solution of the observed system $(\rho,u)$
asymptotically in time. In particular, we are interested in the rate
of convergence of $\| \hat\rho - \rho \|$ and of $\| \hat u - u \|$
towards zero where we will only consider the $L^2$ norm in this paper.

Probably, the simplest observer is defined by Luenberger observer (or
asymptotic observer): as only $u$ is observed, then the feedback term
is added only in the velocity equation ($F_\rho = 0$), and the
feedback term is simply $F_u(\hat\rho,\hat u,u)=k_u(u-\hat u)$, where
$k_u>0$ is a constant. We will study more general observers, with the
aim of correcting also the density equation with the velocity, and in
the process increasing the rate of convergence of the observer towards
the true solution.

In order to study the theoretical behavior of the observer, we first
linearize the system around a steady state of the nonlinear
system~\eqref{eqrhou} and study the convergence of the observer for
the linearized system. We postpone the presentation of numerical
results for the nonlinear case to section~\ref{subsec-nonlin-num}.

\subsection{Linearization around an equilibrium state and linear
  observers} \label{subsec-linear}

We now consider an equilibrium state consisting of constant density
and velocity $(\rho_0,u_0)$, and we linearize equations~\eqref{eqrhou}
around this state:
\begin{equation}
 \begin{array}{l}
   \rho_t + (u_0 \cdot \nabla) \rho + \rho_0 \nabla \cdot u =
      0 \,, \\[0.1cm]
   \rho_0 \left[ u_t + (u_0 \cdot \nabla) u \right] = \mu \Delta u +
     (\lambda + \mu) \nabla (\nabla \cdot u) - \gamma
     \rho_0^{\gamma-1} \nabla \rho \,,\label{eqlrhou1}
 \end{array}
\end{equation}
where only terms linear in $(\rho, u)$ appear. Of course
linearizations around non-constant density and velocity will contain
additional terms such as $(u \cdot \nabla) u_0$ and we will not
consider this case in this paper.

\begin{remark} 
  Note that in dimensions $n \ge 2$, stationary solutions of
  equations~\eqref{eqrhou} could be spatially varying, but in
  one-dimensional case $n=1$, we have the following result:
  \begin{proposition}
    The only stationary solutions of 1D compressible Navier-Stokes
    equations~\eqref{eqrhou} are constants:
    $\rho(t,x) = \rho_0,\ u(t,x)= u_0$.
  \end{proposition}

  \begin{proof}
    Let us write the stationary solution as $\rho(t,x)=\rho_0(x)$ and
    $u(t,x)=u_0(x)$. Then the first of equations~\eqref{eqrhou} in one
    dimension gives $\displaystyle \rho_0(x)={C}/{u_0(x)}$ for a
    constant $C$ determined by the initial conditions. The second of
    equations~\eqref{eqrhou} now reads:
    $$ (C^\gamma u_0^{-\gamma}+Cu_0)_x=\nu_0 u_{0xx} \Leftrightarrow -\gamma
    C^\gamma u_0^{-\gamma-1} u_{0x} + Cu_{0x} = \nu_0 u_{0xx}. $$
    Let us multiply this equation by $u_0$ and integrate over the
    space domain $[0,1]$:
    $$ -\gamma C^\gamma \int_0^1 u_0^{-\gamma} u_{0x} \ dx + C \int_0^1 u_0 u_{0x}
    \ dx = \nu \int_0^1 u_{0xx} u_0\ dx. $$
    The two first integrals are both equal to zero, by considering an
    integration by part, as the boundary conditions are periodic. Then
    $$ 0 = -\nu \int_0^1 (u_{0x})^2\ dx, $$
    meaning that $u_{0x} \equiv 0$.
  \end{proof}
\end{remark}

We can rewrite the solution along the characteristics of the
equation. The transport coefficient for both density and velocity is
$u_0$, hence noting that the derivatives of $\rho(t,x+u_0t)$ are
simply $\nabla \rho(t,x+u_0t)$ and
$\rho_t(t,x+u_0t) + (u_0 \cdot \nabla) \rho(t,x+u_0t)$, and similarly
for derivatives of $u(t,x+u_0t)$, the equations~\eqref{eqlrhou1} at
point $(t,x+u_0t)$ are
\begin{equation}
 \begin{array}{l}
  \rho_t + \rho_0 \nabla \cdot u = 0 \,, \\[0.1cm]
  \rho_0 u_t = \mu \Delta u + (\lambda + \mu) \nabla (\nabla \cdot u)
  - \gamma \rho_0^{\gamma-1} \nabla \rho \,, \label{eql2rhou}
 \end{array}
\end{equation}
which is nothing else than equations~\eqref{eqlrhou1} with $u_0=0$. So
we can now assume that $u_0=0$ without any loss of generality. Thus we
will work with this linearized Navier-Stokes system with initial
conditions $\rho(0,x) = \rho_I(x)$ and $u(0,x)=u_I(x)$, and periodic
boundary conditions as in equation~\eqref{eqbc}.

As stated above, we assume that the initial conditions $\rho_I(x)$ and
$u_I(x)$ are unknown, but we have observations of the solution
$u(t,x)$ of the linearized equations~\eqref{eql2rhou}. Again, the goal
is to build an observer $(\hat\rho,\hat u)$ for this system, in such a
way that the observer $(\hat\rho,\hat u)$ converges towards the
solution $(\rho,u)$.  We will use the same observer as introduced in
equations~\eqref{eqobsrhou}-\eqref{feedbackrhou}, except that the left
hand side is now linear, just as in equations~\eqref{eql2rhou}:
\begin{equation}
 \begin{array}{l}
   \hat\rho_t + \rho_0 \nabla \cdot \hat u = \varphi_\rho\ast D_\rho
    \left[ (u-\hat u) \mathbbm{1}_\Omega \right] \,, \\[0.1cm]
   \rho_0 \hat u_t = \mu \Delta \hat u + (\lambda + \mu) \nabla
    (\nabla \cdot \hat u) - \gamma \rho_0^{\gamma-1} \nabla \hat \rho
    + \varphi_u\ast D_u\left[(u-\hat u) \mathbbm{1}_\Omega \right]
    \,, \label{eqobslrhou}
 \end{array}
\end{equation}
with periodic boundary conditions, and initial conditions different
from the true initial conditions:
$\hat\rho(0,x) = \hat\rho_I(x) \ne \rho_I(x)$ and
$\hat u(0,x) = \hat u_I(x) \ne u_I(x)$.

We will now present the main theoretical result about this observer in
the linear case, stating that we can choose the feedback terms in
order to guarantee any specified rate of convergence for arbitrarily
large but finitely many number of Fourier modes. Since we are dealing
with linear PDE, we will use Fourier series representation as our main
tool.

\section{Theoretical study of an observer for n-dimensional linear
  Navier-Stokes system} \label{sec-theory}

For the case of velocity observations over the full domain
($\Omega=[0,1]^n$), subtracting the observer
equations~\eqref{eqobslrhou} from the state
equations~\eqref{eql2rhou}, we get the following equations for the
errors $r = \hat\rho - \rho$ and $v = \hat u - u$:
\begin{equation}
 \begin{array}{l}
  r_t + \rho_0 \nabla\cdot v = -\varphi_\rho \ast D_\rho v \,,
   \\[0.1cm] \rho_0 v_t = \mu \Delta v +
   (\lambda+\mu)\nabla(\nabla\cdot v)-\gamma \rho_0^{\gamma-1}\nabla r
   - \varphi_u\ast D_u v \,, \label{erroreq} 
 \end{array}
\end{equation}
which are exactly identical to equations~\eqref{eqobslrhou}, with
$(r,v)$ replacing $(\hat\rho, \hat u)$ and with $u = 0$.

\subsection{Damped wave equation formulation} \label{subsec-damped}

We now eliminate the density in the velocity equation, in order to get
an equation for the velocity alone. As we will see, this will allow us
to more easily define the observers, and more particularly the
differential operators $D_\rho$ and $D_u$. Starting from
equations~\eqref{erroreq}, we take the space derivative of the density
equation, and the time derivative of the velocity equations:
$$ \begin{array}{l} \nabla r_{t}+\rho_0 \nabla(\nabla\cdot v) =
     -\varphi_\rho\ast  \nabla (D_\rho v),\\ \rho_0 v_{tt}=\mu\Delta
     v_t + (\lambda+\mu)\nabla(\nabla\cdot
     v_t)-\gamma\rho_0^{\gamma-1}\nabla r_t - \varphi_u \ast D_u v_t
     \,. \end{array} $$
We replace $\nabla r_{t}$ in the second equation by its expression
given by the first equation, and we obtain:
\begin{equation}
  \rho_0 v_{tt}=\mu\Delta v_t+(\lambda+\mu)\nabla(\nabla\cdot v_t)+\gamma\rho_0^\gamma \nabla(\nabla\cdot v)+\gamma\rho_0^{\gamma-1}\varphi_\rho\ast\nabla(D_\rho v)-\varphi_u \ast D_u v_t\,.
\label{eqwave}
\end{equation}

Equation~\eqref{eqwave} is a damped wave equation, with two forcing
terms coming from the observer feedbacks. In this case, the goal is to
make the difference $v = \hat{u} - u$ between the observer $\hat{u}$
and the observed velocity $u$ converge towards $0$.  We want to choose
$D_\rho$ and $D_u$ which respect the symmetries of the original
system, but of course the choice is not unique. Since the observations
may be noisy, we also look for the operators with minimum number of
derivatives. With this in mind, we now assume that the following
differential operators are used:
\begin{equation}\label{diffop}
  D_\rho(f) = \rho_0 \nabla\cdot f,\quad D_u(f) = \rho_0 f,
\end{equation}
which means that, up to constants, we want the velocity equation to be
controlled by the velocity, and the density equation by the divergence
of the velocity. In the damped wave equation~\eqref{eqwave},

Using this choice, equation~\eqref{eqwave} rewrites as follows:
\begin{equation}
  \rho_0 v_{tt}=\mu\Delta v_t+(\lambda+\mu)\nabla(\nabla\cdot v_t)+\gamma\rho_0^\gamma \nabla(\nabla\cdot v)+\gamma\rho_0^{\gamma}\varphi_\rho\ast\nabla(\nabla\cdot v)-\rho_0\varphi_u \ast v_t\,,
\label{eqwave2}
\end{equation}
from where we can see that the last term (with $\varphi_u$) controls
the first two terms on the right while the fourth term (with
$\varphi_\rho$) controls the third term.

Using equations~\eqref{diffop}, the observer error
equations~\eqref{erroreq} now read:
\begin{equation}
 \begin{array}{l}
  r_t + \rho_0 \nabla\cdot v = -\rho_0 \varphi_\rho \ast \nabla\cdot v \,, \\[0.1cm]
  \rho_0 v_t = \mu \Delta v + (\lambda+\mu)\nabla(\nabla\cdot v)-\gamma \rho_0^{\gamma-1}\nabla r - \rho_0 \varphi_u\ast v. \label{erroreq2}
 \end{array}
\end{equation}
We will refer to this system as the full state observer, as it
corrects both velocity and density equations.

\subsection{Fourier transform} \label{subsec-fourier}

We first note that equation~\eqref{erroreq2} defines the solution up
to a constant. This, if the mean value $m$ of $r(t,x)$ is not equal to
$0$, then the solution can at best converge towards the constant $m$
instead of $0$. Hence we assume the following.
\begin{assumption}\label{assumption}
  We suppose that the mean value of $r(t,x)$ is equal to $0$, i.e. the
  mean values of the initial conditions of $\hat\rho$ and $\rho$ are
  the same.
\end{assumption}

Of course, this is a strong assumption as the initial state could have
bias. But it is quite often assumed, at least in the data assimilation
community, that bias on observations could be removed, for instance
with reanalysis methods, and that noises are zero-mean. Moreover, the
linear system~\eqref{eqlrhou1} defines the velocity up to a constant,
and the observer that we propose is consistant with this. We also note
that we will be able to control the mean value of the velocity, as
shown below.

As we are on the domain $[0,1]^n$ with periodic boundary conditions
(see equation~\eqref{eqbc}), we can consider the following Fourier
decomposition of the velocity and density:
\begin{equation}
  v(t,x)=\sum_{k} a_k(t) e^{i2\pi k\cdot x}, \qquad
  r(t,x)=\sum_{k} b_k(t) e^{i2\pi k\cdot x},
\label{eqfourier}
\end{equation}
where $a_k(t)\in\mathbb{R}^n$ and $b_k(t)\in\mathbb{R}$ are the time
dependent Fourier coefficients,
$$a_k(t)=
\left(\begin{array}{c}
a_{1k}(t)\\ \vdots \\ a_{nk}(t)
\end{array}\right), \quad \textrm{and} \quad 
k=\left(\begin{array}{c}k_1\\ \vdots \\k_n\end{array}\right)\in\mathbb{Z}^n.$$

Note that we assume $b_0=0$, from assumption~\ref{assumption}.  We
denote by $\varphi_{\rho k}$ and $\varphi_{u k}$ the Fourier
coefficients of the (time independent) functions $\varphi_\rho(x)$ and
$\varphi_u(x)$ respectively.

Substituting~\eqref{eqfourier} in~\eqref{erroreq2}, we obtain the following equations for mode $k$: for all $1\le i\le n$,
\begin{equation}
 \rho_0 a_{ik}'(t)=-(\mu 4\pi^2|k|^2+\rho_0\varphi_{uk})a_{ik}(t)-(\lambda+\mu)4\pi^2\left(\sum_{j=1}^na_{jk}(t)k_j\right)k_i -\gamma\rho_0^{\gamma-1} i2\pi b_k(t)k_i, \label{eqmode1a}
\end{equation}
where $|k|^2=\sum_{j=1}^n k_j^2$, and
\begin{equation}
 b_k'(t)=-\rho_0 i2\pi (1+\varphi_{\rho k}) \sum_{j=1}^n a_{jk}(t)k_j.
\label{eqmode1b}
\end{equation}

Introducing $y_k(t)=\left(\begin{array}{c}b_k(t)\\a_{1k}(t)\\ \vdots \\a_{nk}(t)\end{array}\right)\in\mathbb{R}^{n+1}$, the system~\eqref{eqmode1a}-\eqref{eqmode1b} rewrites as 
\begin{equation}\label{eqy}
y'_k(t)=M_ky_k(t),
\end{equation}
where $M_k$ is the following $(n+1)\times (n+1)$ complex matrix:
\begin{equation}
M_k=\left( \begin{array}{cc}
0 & c_4 (1+\varphi_{\rho k})k^T \\
- c_3 k & -(c_1|k|^2+\varphi_{uk})I_n-c_2 k\otimes k
\end{array}\right)
\label{eq:nk} \end{equation}
where $I_n$ is the $n\times n$ identity matrix, $0$ is the scalar zero, and
\begin{equation}\label{c14eq}
c_1=\frac{\mu 4\pi^2}{\rho_0},\quad c_2=\frac{(\lambda+\mu)4\pi^2}{\rho_0}, \quad c_3=\gamma \rho_0^{\gamma -2}i 2 \pi, \quad c_4=-\rho_0 i 2 \pi \,.
\end{equation}

\begin{remark}
Note that the Fourier analysis can also be done on the damped wave equation~\eqref{eqwave}. Substituting~\eqref{eqfourier} in~\eqref{eqwave}, we obtain the
following equation for mode $k$: for all $1\le i\le n$,
\begin{eqnarray}
 \rho_0 a_{ik}''(t)=-(\mu 4\pi^2|k|^2+\rho_0\varphi_{uk})a_{ik}'(t)&&-(\lambda+\mu)4\pi^2\left(\sum_{j=1}^na_{jk}'(t)k_j\right)k_i\nonumber\\
&&-\gamma\rho_0^\gamma(1+\varphi_{\rho k}) 4\pi^2\left(\sum_{j=1}^n a_{jk}(t)k_j\right)k_i.
\label{eqmode}
\end{eqnarray}
Let $z_k(t)=\left(\begin{array}{c}a_{1k}(t)\\ \vdots \\ a_{nk}(t)\\a_{1k}'(t)\\ \vdots \\ a_{nk}'(t)\end{array}\right)\in\mathbb{R}^{2n}$, then equation~\eqref{eqmode} rewrites as $z_k'(t)=N_k z_k(t)$, where $N_k$ is the following $2n\times 2n$ matrix:
\begin{equation}
N_k=\left( \begin{array}{cc}
0_n & I_n \\
-c_5(1+\varphi_{\rho k})k\otimes k & -(c_1|k|^2+\varphi_{uk})I_n-c_2 k\otimes k
\end{array}\right)
\label{eq:mk} \end{equation}
where $0_n$ is the $n\times n$ zero matrix, and where we recall de 
values of $c_1$ and $c_2$ and introduce $c_5$:
\begin{equation}
c_1=\frac{\mu 4\pi^2}{\rho_0}, \quad
c_2=\frac{(\lambda+\mu)4\pi^2}{\rho_0}, \quad
c_5 =c_3\times 
c_4=\gamma\rho_0^{\gamma -1}4\pi^2 \, .
\end{equation}
\end{remark}

\subsection{Spectral analysis}\label{subsec-spectral}

We now look at the eigenvalues and eigenvectors of matrix $M_k$. For $k\ne 0$, Let $\{k_i^\perp,1\le i\le n-1\}$ be a basis of the set of vectors orthogonal to $k$ in $\mathbb{R}^n$.

\begin{proposition}
\begin{itemize}
\item For $k\ne 0$, the eigenvalues of matrix $M_k$ are:
\begin{equation} \label{eqeigN}
 \begin{array}{l}
  \lambda_{dk}=-(c_1|k|^2+\varphi_{uk}),\\[0.1cm]
  \displaystyle \lambda_{\pm k}=\frac{-\left((c_1+c_2)|k|^2+\varphi_{uk}\right)\pm\sqrt{\left((c_1+c_2)|k|^2+\varphi_{uk}\right)^2-4c_5(1+\varphi_{\rho k})|k|^2}}{2},
 \end{array}
\end{equation}
with
\begin{equation}\label{cvalues}
c_1=\frac{\mu 4\pi^2}{\rho_0}, \quad
c_2=\frac{(\lambda+\mu)4\pi^2}{\rho_0}, \quad
c_5 =\gamma\rho_0^{\gamma -1}4\pi^2 \, .
\end{equation}
$\lambda_{dk}$ has a multiplicity $n-1$, and the associated eigenvectors are the following $n-1$ vectors in $\mathbb{R}^{n+1}$: $\displaystyle \left(\begin{array}{c}0\\k_i^\perp\end{array}\right)$.\\
$\lambda_{+k}$ and $\lambda_{-k}$ both have a multiplicity $1$, and the associated eigenvectors are respectively: $\displaystyle \left(\begin{array}{c}c_4(1+\varphi_{\rho k})|k|^2\\ \lambda_{+k}k\end{array}\right)\quad \textrm{and} \quad \left(\begin{array}{c}c_4(1+\varphi_{\rho k})|k|^2\\ \lambda_{-k}k\end{array}\right)$.
\item For $k=0$, the eigenvalues of $M_0$ are:
\begin{equation}
\lambda_{d0}=-\varphi_{u0}, \quad \lambda_{+0}=0.
\end{equation}
$\lambda_{d0}$ now has multiplicity $n$, with the associated eigenvectors: $\displaystyle\left(\begin{array}{c}0\\e_i\end{array}\right)$, $\{e_i,1\le i\le n\}$ being a basis of $\mathbb{R}^n$, and $\lambda_{+0}$ has multiplicity $1$, with the associated eigenvector $\displaystyle\left(\begin{array}{c}1\\0_{\mathbb{R}^n}\end{array}\right)$.
\end{itemize}
\end{proposition}

The proof is straightforward and simply consists in multiplying these eigenvectors by $M_k$. Note that in the case $k=0$, the eigenvalue $\lambda_{+0}$ cannot be controlled, and the associated eigenvector corresponds to constant density solutions. So, from assumption~\ref{assumption}, the projection of $r(x,t), v(x,t)$ along this eigenvector will be zero.

For $k\ne 0$, the discriminant of the quadratic equation for the eigenvalues $\lambda_{\pm k}$ is
\begin{equation}\label{eqdelta}
\Delta_k=\left((c_1+c_2)|k|^2+\varphi_{uk}\right)^2-4c_5(1+\varphi_{\rho k})|k|^2
\end{equation}
and we use the notation $\sqrt{\Delta_k}=i\sqrt{-\Delta_k}$ if the discriminant
is negative.

The eigenvectors associated to $\lambda_{dk}$ correspond to non constant divergence free velocity solutions, for which the equation simply consists in a diffusion equation (plus the feedback term), the diffusion coefficient being $\lambda_{dk}$. Also note that the first component of the two eigenvectors associated to $\lambda_{\pm k}$ is complex, since $c_4 = -\rho_0 i2\pi$.

We will now relate the eigenvalues of the matrix $M_k$ to the \emph{decay rate} of the solutions $y_k(t)$ of equation~\eqref{eqy}.
\begin{definition}
$d_k>0$ is a decay rate of $y_k(t)$ if there exists a constant $c>0$ such that $\|y_k(t)\|\le ce^{-d_kt}$ for any $t\ge 0$.
\end{definition}

Of course, we are usually interested in the largest decay rate. In our particular situation, the eigenvalues given by equations~\eqref{eqeigN} are all negative if $\varphi_{uk}\ge 0$ and $\varphi_{\rho k}\ge 0$. In such a case, it is straightforward to see that the largest decay rate is given by:
\begin{equation}
d_k=-\max_{\lambda\in Sp(M_k)}\left\{ \Re(\lambda)\right\}
\end{equation}
because there exists $c>0$ such that $\displaystyle \|y_k(t)\|\le c e^{\max\{\Re(\lambda)\}t}$. Thus, the largest decay rate can be controlled by $\varphi_{uk}$ and $\varphi_{\rho k}$ as we will see in next section.

\begin{remark}
For $k\ne 0$, the eigenvalues of matrix $N_k$ are the following:
\begin{equation} \lambda_{0k}=0, \quad \lambda_{dk}, \quad \lambda_{\pm k},\label{eqeigM}
\end{equation}
where $\lambda_{dk}$ and $\lambda_{\pm k}$ are given by~\eqref{eqeigN}.

The new eigenvalue $\lambda_{0k}$ has a multiplicity $n-1$, and $\lambda_{dk}$ still has a multiplicity $n-1$. The $n-1$ associated eigenvectors are:
 $\displaystyle \left(\begin{array}{c}k_i^\perp \\ \lambda k_i^\perp \end{array}\right)$, with $\lambda=\lambda_{0k}$ or $\lambda_{dk}$. 

Note that the solutions corresponding to $\lambda_{0k}$ are solutions of the damped wave equation~\eqref{eqwave}, but not solutions of the original system~\eqref{erroreq}.

These eigenvectors correspond to divergence free solutions, either constant in time (for $\lambda=\lambda_{0k}$) or not (for $\lambda=\lambda_{dk}$).

Finally, $\lambda_{+k}$ and $\lambda_{-k}$ both have a multiplicity $1$, and the eigenvectors are respectively: $\displaystyle \left(\begin{array}{c}k\\ \lambda_{+k}k\end{array}\right)\quad \textrm{and} \quad \left(\begin{array}{c}k\\ \lambda_{-k}k\end{array}\right)$.

Note that the eigenvalue $\lambda_{0k}$ is artificial, and only appears because we took the time derivative of the velocity equation in order to eliminate the density. This is the only main difference between $M_k$ and $N_k$ spectral analysis.
\end{remark}

\subsection{Main result} \label{subsec-result}

We now give the main result in this framework: 
\begin{theorem} \label{mainprop} We assume that
  assumption~\ref{assumption} holds. For any $d>0$, for any $K>0$, one
  can find $\varphi_{\rho}(x)$ and $\varphi_{u}(x)$ such that the
  maximal decay rate of the errors $r(t,x)$ and $v(t,x)$, solutions
  of~\eqref{erroreq2}, towards $0$ is at least $d$ for any Fourier
  mode $k$ such that $|k|\le K$. The following values can be chosen,
  with $c_1$, $c_2$ and $c_5$ given by~\eqref{cvalues}:
  \begin{eqnarray}
    \varphi_{u0} &=& d \,, \quad \varphi_{\rho 0} = 0 \,, \label{varphiu0}\\
    \varphi_{uk} &=& \max\{0;d-c_1|k|^2;2d-(c_1+c_2)|k|^2)\}\label{varphiu} \,, \quad 0 < |k| \le K \,,\\
    \varphi_{\rho k} &=& \max\left\{0;\frac{\left((c_1+c_2)|k|^2+\varphi_{uk}\right)^2}{4c_5|k|^2}-1\right\}\,, \quad  0 < |k| \le K \,,\label{varphirho} \\
    \varphi_{uk} &=& 0 \,, \quad \varphi_{\rho k} = 0 \,, \quad |k| > K \,.
  \end{eqnarray}
\end{theorem}

In other words, for any specified decay rate, we can find convolution
kernels $\varphi_\rho$ and $\varphi_u$ such that the observer $\hat u$
converges towards $u$ at this specified rate up to any Fourier
mode. Indeed, we can choose appropriate Fourier coefficients of these
kernels for any mode $k \in \mathbb{Z}^n$, but as we will see, their
expression does not ensure the convergence of the Fourier series and
we need to truncate the series.

\begin{proof}
  Let $d>0$, $K>0$, and let $k$ a Fourier mode such that $0<|k|\le
  K$.

  The case $k=0$ is trivial: the only eigenvalue is
  $\lambda_{d0} = -\varphi_{u0}$, leading to the appropriate decay
  rate for the velocity if we choose equation~\eqref{varphiu0}. There
  is no correction on the density but the Fourier mode $b_0$ is $0$
  from assumption~\ref{assumption}.
  
  We can see that if we use~\eqref{varphiu}, then this choice ensures
  non-negativity of all Fourier coefficients of the convolution
  kernels, and from~\eqref{eqeigN}, we see that $-\lambda_{dk}$ is at
  least $d$. For small modes, the diffusion process is not large
  enough to ensure the decay rate, so we need to add the feedback
  term. For larger modes, diffusion will be enough, and there is no
  need to add the feedback (but one can still add a feedback term, and
  the decay rate will become larger than the specified rate for such
  modes). Note that even if we drop the $\max\{0;.\}$
  in~\eqref{varphiu}, it still ensures that the decay rate will be at
  least $d$, but $\varphi_{uk}$ may become negative for large modes,
  and this can lead to numerical instabilities.

  Then, recalling equation~\eqref{eqdelta}, we set $\varphi_{\rho k}$
  as in~\eqref{varphirho}. This choice of $\varphi_{\rho k}$ ensures
  that $\Delta_k \le 0$, and the decay rate corresponding to
  eigenvalues $\lambda_{\pm k}$ is then exactly
  $\displaystyle \frac{(c_1+c_2)|k|^2+\varphi_{uk}}{2}$, which is
  always larger than (or equal to) $d$, because of the choice of
  $\varphi_{uk}$ in equation~\eqref{varphiu}.

  The $\max\{0,.\}$ in equation~\eqref{varphirho} is also optional,
  since even without the $\max$, the discriminant $\Delta_k$ will be
  negative, but the non-negativity of $\varphi_{\rho k}$ will avoid
  numerical instabilities.

  For large modes, namely $|k|>\sqrt{\frac{d}{c_1}}$ and
  $|k|>\sqrt{\frac{2d}{c_1+c_2}}$, $\varphi_{uk}=0$, and then there is
  no issue in considering the inverse Fourier transform and defining a
  convolution kernel $\varphi_u(x)$ with these coefficients.

  With the above choices, only finitely many Fourier modes of
  $\varphi_u(x)$ and $\varphi_\rho(x)$ are nonzero and hence the
  convolution kernels exist, which proves the result.
\end{proof}

Note that we could use the same definitions from
equations~\eqref{varphiu}-\eqref{varphirho} for larger modes
($|k|>K$), but then the Fourier series for $\varphi_\rho(x)$ does not
converge, as $\varphi_{\rho k}=\mathcal{O}(|k|^2)$. So we need to
truncate the series, and set $\varphi_{\rho k}=0$ for $|k|>K$.

We also note that for large modes for which both $\varphi_{\rho k}$
and $\varphi_{uk}$ are set to 0, the observer is simply a solution of
the equation without any forcing term, and the largest eigenvalue is
$\lambda_{+k}$ (both $\lambda_{dk}$ and $\lambda_{-k}\to -\infty$ when
$|k|\to +\infty$). From~\eqref{eqeigN}, the decay rate is then equal
to
$\displaystyle\frac{c_5}{c_1+c_2}=\frac{\gamma\rho_0^\gamma}{\lambda+2\mu}$
asymptotically for $|k| \to +\infty$.

\subsection{Observers with density
  observations} \label{subsec:rho-obs}

We have only considered observations of the velocity field so far in
this paper.  In this subsection, we will construct the observers for
the case when the density is observed instead of velocity.

Since observations of the scalar field (density) naturally gives less
information than those of a vector field (velocity) in any dimension
greater than one, we will see that, as expected, we will not be able
to control all the modes in the observer. In particular, we will be
unable to modify the convergence rates for the eigenvectors
corresponding to the diffusion eigenvalues but we will be able to get
arbitrarily high convergence rates for the eigenvectors corresponding
to eigenvalues $\lambda_\pm$. Since the details of the calculations
are very similar to those presented above, we will only present a
sketch of the results.

When we have density observations $\rho(t,x)$ instead of velocity
$u(t,x)$, the observer is exactly as in the observer
equations~\eqref{eqobslrhou} with $u$ being replaced by $\rho$. Then
the equations for errors $r = \hat\rho - \rho$ and $v = \hat u - u$
will be as follows. These are the same as equations~\eqref{erroreq} with $r$ replacing $v$ in the feedback terms.
\begin{equation}
 \begin{array}{l}
   r_t + \rho_0 \nabla\cdot v = -\psi_\rho \ast \bar{D}_\rho r\,, \\[0.1cm]
   \rho_0 v_t = \mu \Delta v + (\lambda+\mu)\nabla(\nabla\cdot v)-\gamma \rho_0^{\gamma-1}\nabla r - \psi_u\ast \bar{D}_u r \,. \label{erroreqrho}
 \end{array}
\end{equation}
We see that because of the symmetries of the
equations, the natural choice of the differential operators
$\bar{D}_\rho$ and $\bar{D}_u$ with minimal derivatives is
\begin{equation}\label{barD}
  \bar{D}_\rho(f) = f \quad \textrm{and} \quad \bar{D}_u(f) = \gamma \rho_0^{\gamma-1} \nabla f\,.
\end{equation}

\begin{assumption}\label{assumption2}
We suppose that the mean value of $v(t,x)$ is equal to $0$, i.e. the mean values of the initial conditions of $\hat u$ and $u$ are the same.
\end{assumption}

Again, this is a strong assumption similar to assumption~\ref{assumption}: contrary to the density mean, the velocity mean cannot be corrected from~\eqref{erroreqrho}-\eqref{barD}.

\begin{remark}
As in the case of velocity observations, we can eliminate $v$ in order to obtain an equation for the density alone. This can help understanding the choice of the differential operators in the feedback terms. From the first equation of~\eqref{erroreqrho}, we can extract 
$\rho_0 \nabla\cdot v=-r_t-\psi_\rho\ast \bar{D}_\rho r$ and then we need
to take the divergence of the second equation:
$$ \rho_0 \nabla\cdot v_t = \mu \nabla\cdot \Delta v + (\lambda+\mu)\nabla\cdot (\nabla(\nabla\cdot v))-\gamma\rho_0^{\gamma-1}\nabla\cdot(\nabla r)-\psi_u\ast \nabla\cdot (\bar D_u r).$$
By choosing $\bar{D}_\rho$ and $\bar{D}_u$ as in~\eqref{barD}, we get the damped wave equation for $r(x,t)$
\begin{equation}
  \rho_0 r_{tt} = (\lambda+2\mu) \Delta r_t + \gamma \rho_0^\gamma \Delta r + (\lambda+2\mu)\psi_\rho\ast\Delta r + \gamma\rho_0^\gamma \psi_u\ast\Delta r -\rho_0\psi_\rho\ast r_t\,,
\label{eqsym}
\end{equation}
which again shows that the choices of $\bar{D}_u$ and $\bar{D}_\rho$
are quite natural to control the terms in the above equation.

We can then see that if we are in one dimension $n=1$, equations~\eqref{eqwave2} and~\eqref{eqsym} are identical with $\varphi_u = \psi_\rho$, and
$\varphi_\rho = \psi_u + ((2\mu+\lambda)\psi_\rho) / (\gamma\rho_0^\gamma)$. Thus in one dimension, observations of the density and those of velocity give the same results for rates of convergence of the observer solution. 
\end{remark}

Introducing the Fourier decomposition for $(r,v)$ exactly like in
equations~\eqref{eqfourier} and the vector
$y_k$ of the Fourier coefficients, we will get a system
$y'_k(t)=\bar{M}_ky_k(t)$, where $\bar{M}_k$ is the following
$(n+1)\times (n+1)$ complex matrix:
\begin{equation}
\bar{M}_k=\left( \begin{array}{cc}
-\psi_{\rho k} & c_4 k^T \\
- c_3(1 + \psi_{uk}) k & -c_1|k|^2 I_n-c_2 k\otimes k
\end{array}\right)
\label{eq:nkbar} \end{equation}
We recall the values of the constants:
\begin{equation}\label{cvalues2}
c_1=\frac{\mu 4\pi^2}{\rho_0},\ c_2=\frac{(\lambda+\mu)4\pi^2}{\rho_0}, \ c_3=\gamma \rho_0^{\gamma -2}i 2 \pi, \ c_4=-\rho_0 i 2 \pi, \ c_5=c_3\times c_4 = \gamma \rho_0^{\gamma -1}4\pi^2.
\end{equation}
We also recall that $\{k_i^\perp,1\le i\le n-1\}$ is a basis of the set of vectors orthogonal to $k$ in $\mathbb{R}^n$ if $k\ne 0$. We have then the following result:

\begin{proposition}
\begin{itemize}
\item For $k\ne 0$, the eigenvalues of matrix $\bar{M}_k$ are:
\begin{equation}
 \begin{array}{l}
  \bar\lambda_{dk} = -c_1 |k|^2, \\[0.1cm]
  \displaystyle\bar\lambda_{\pm k}=\frac{-\left((c_1+c_2)|k|^2+\psi_{\rho k}\right)\pm\sqrt{\left((c_1+c_2)|k|^2-\psi_{\rho k}\right)^2-4c_5(1+\psi_{uk})|k|^2}}{2} \,.
 \end{array}
\label{eqeigrho}
\end{equation}
The eigenvalue $\bar\lambda_{dk}$ has a multiplicity $n-1$, with eigenvectors $\left(\begin{array}{c}0\\k_i^\perp\end{array}\right)$. The eigenvalues $\bar\lambda_{\pm k}$ both have a multiplicity $1$, with eigenvectors $\left(\begin{array}{c}c_4|k|^2\\(\lambda_{\pm k}+\psi_{\rho k})k\end{array}\right)$.
\item For $k=0$, the eigenvalues of matrix $\bar{M}_0$ are:
\begin{equation}
\bar\lambda_{d0}=0, \quad \bar\lambda_{-0}=-\psi_{\rho 0}.
\end{equation}
$\bar\lambda_{d0}$ now has multiplicity $n$, with associated eigenvectors: $\displaystyle\left(\begin{array}{c}0\\e_i\end{array}\right)$, $\{e_i,1\le i\le n\}$ being a basis of $\mathbb{R}^n$, and $\bar\lambda_{-0}$ has multiplicity $1$, with the associated eigenvector $\displaystyle\left(\begin{array}{c}1\\0_{\mathbb{R}^n}\end{array}\right)$.
\end{itemize}
\end{proposition}

The proof of this result is straightforward, one simply has to multiply these eigenvectors by $\bar{M}_k$.

Note that for $k=0$, $\bar\lambda_{d0}$ cannot be controlled, but the associated eigenvectors correspond to constant velocity solutions. So from assumption~\ref{assumption2}, this case does not appear. But in this case $\bar\lambda_{-0}$ can be controlled, and it corresponds to constant density solutions.

For $k\ne 0$, the discriminant is
\begin{equation}\label{eqdeltarho}
\bar\Delta_k=\left((c_1+c_2)|k|^2-\psi_{\rho k}\right)^2-4c_5(1+\psi_{uk})|k|^2
\end{equation}
and we still use the notation $\sqrt{\bar\Delta_k}=i\sqrt{-\bar\Delta_k}$ if $\bar\Delta_k < 0$.

Note that $\bar\lambda_{dk}$ is independent of $\psi_\rho$ and $\psi_u$ and
thus cannot be controlled by any choice of $\psi_\rho$ or $\psi_u$. Exactly in parallel with the previous case (velocity observations) and the choice presented in equations~\eqref{varphiu}-\eqref{varphirho}, we have the following result:
\begin{theorem}\label{proprho}
  We suppose that assumption~\ref{assumption2} holds. For any $d>0$,
  for any $K>0$, one can find $\psi_{\rho}(x)$ and $\psi_{u}(x)$ such
  that the maximal decay rate of the errors $r(t,x)$ and $v(t,x)$,
  solutions of~\eqref{erroreqrho}, towards $0$ is at least
  $\min\{d,c_1|k|^2\}$ for any Fourier mode $k$ such that $|k|\le
  K$.
  The following values can be chosen, with $c_1$, $c_2$ and $c_5$
  given by~\eqref{cvalues2}:
  \begin{eqnarray}
    \psi_{\rho 0} &=& d \,, \quad \psi_{u0} = 0 \,, \label{psirho0}\\
    \psi_{\rho k} &=& \max\{0;2d-(c_1+c_2)|k|^2)\} \,, \quad 0 < |k| \le K \,,\label{psirho}\\
    \psi_{uk} &=& \max\left\{0;\frac{\left((c_1+c_2)|k|^2-\psi_{\rho k}\right)^2}{4c_5|k|^2}-1\right\} \,, \quad 0 < |k| \le K \,, \label{psiu} \\
    \psi_{uk} &=& 0 \,, \quad \psi_{\rho k} = 0 \,, \quad |k| > K \,.
  \end{eqnarray}
\end{theorem}
The $\max(.,0)$ in~\eqref{psiu} is there to ensure positiveness of all Fourier coefficients of $\psi_u$. It could be relaxed as it is only useful for small modes (the other term grows as $|k|^2$).
\begin{proof}
The choice of $\psi_{uk}$ ensures that $\bar\Delta_k \le 0$ and the real part of
$\bar\lambda_{\pm k}$ is then $\displaystyle-\frac{\left((c_1+c_2)|k|^2+\psi{\rho k}\right)}{2}$. Then, using~\eqref{psirho}, the decay rate is then at least $d$ for all $k$ with $|k| \le K$.
\end{proof}

We note that the modes that we cannot control to have a pre-specified
decay rate are the incompressible ones since the eigenvector
$\left(\begin{array}{c}0\\k_i^\perp\end{array}\right)$ satisfies
$k \cdot v_k = 0$, corresponding to $\nabla \cdot v = 0$. This is
expected since the density observations cannot give any information
about the ``constant density'' incompressible flow. 

\subsection{Remarks and comparison with the nudging feedback} \label{subsec-remark}

We now compare the result presented in the previous sections with what
can be done with the standard nudging observer, for which only the observed variable is corrected. We will
use $a=(c_1+c_2)|k|^2$, $b=2\sqrt{c_5}|k|$  (see~\eqref{cvalues2}) in the discussion in the
rest of this section.

\paragraph{Nudging with velocity observations}
In this case, there is no feedback on the density: $\varphi_\rho=0$, and the differential operator on the velocity $D_u$ is simply the identity (or more physically $D_u(f)=\rho_0 f$ as in~\eqref{diffop}). Then equations~\eqref{erroreq} become 
\begin{equation}
 \begin{array}{l}
  r_t + \rho_0 \nabla\cdot v = 0 \,, \\[0.1cm]
  \rho_0 v_t = \mu \Delta v + (\lambda+\mu)\nabla(\nabla\cdot v)-\gamma \rho_0^{\gamma-1}\nabla r - \rho_0 \varphi_u\ast v \,. \label{erroreqnudg1}
 \end{array}
\end{equation}
Let us fix
$x=\varphi_{uk}$ for this paragraph. We can write the eigenvalues
$\lambda_{\pm}$ from equation~\eqref{eqeigN} as
\begin{align*}
  \lambda_{\pm} = -\frac{1}{2} \left[(a+x) \pm \sqrt{(a+x)^2 - b^2} \right].
\end{align*}
From this expression, we can find optimal decay rates as follows.
\begin{itemize}
\item If $a \le b$, i.e, for $|k| \le \frac{2\sqrt{c_5}}{c_1+c+2}$, it
  is easy to see that the real part of $-\lambda_{\pm}$ is maximized
  for $x = b-a$, and for this choice the decay rate is
  $b/2 = \sqrt{c_5}|k|$ which is greater than the decay rate without
  nudging ($x=0$) but it cannot be made arbitrarily large as in the full state observer (see theorem~\ref{mainprop}).
\item If $a > b$, i.e., for $|k| > \frac{2\sqrt{c_5}}{c_1+c+2}$, we
  again see that the real part of $-\lambda_{\pm}$ is maximized for
  $x=0$ and this simply gives the decay rate without nudging. Again,
  the full state observer presented above is better then simple nudging since we can obtain
  arbitrarily large decay rates up to $|k| < K$ for any fixed $K$.
\end{itemize}

\paragraph{Nudging with density observations} 
In this case, there is no feedback term on the velocity: $\psi_u=0$, and the differential operator on the density $\bar{D}_\rho$ is still the identity, as in~\eqref{barD}. The system~\eqref{erroreqrho} now becomes
\begin{equation}
 \begin{array}{l}
   r_t + \rho_0 \nabla\cdot v = -\psi_\rho \ast  r\,, \\[0.1cm]
   \rho_0 v_t = \mu \Delta v + (\lambda+\mu)\nabla(\nabla\cdot v)-\gamma \rho_0^{\gamma-1}\nabla r \,. \label{erroreqrho2}
 \end{array}
\end{equation}
Let us fix
$x=\psi_{\rho k}$ for this paragraph. We can write the eigenvalues
$\bar\lambda_{\pm}$ from equation~\eqref{eqeigrho} as
\begin{align*}
  \bar\lambda_{\pm} = -\frac{1}{2} \left[(a+x) \pm \sqrt{(a-x)^2 - b^2} \right].
\end{align*}
From this expression, we can find optimal decay rates by maximizing
the real part of the above expression as a function of $x$ in exactly
the same manner as the discussion in the previous paragraph.  In this
case, it is easy to see that the real part of $-\lambda_{\pm}$ is
maximized for $x = b+a$, and for this choice the decay rate is
$\min\{a + b/2, c_1|k|^2\}$ which is greater than the decay rate
without nudging ($x=0$) but the first factor in the $\min$ cannot be
made arbitrarily large, and the decay rate may be smaller than the one for
the full state observer (see theorem~\ref{proprho}).

Thus we see that in both cases of either velocity or density
observations, the full state observers previously presented 
perform better than simple nudging, even with an optimal choice of the nudging
convolution kernel.

\section{Additional study of the observer in dimension 1} \label{sec-theory-1d}

In this section, we study the same observer as in the above
discussion, but in one dimension. There are a few simplification which
we point out below. We extend the study of the above observer to two
special cases in one dimension. In subsection~\ref{subsec-forcing}, we
discuss the case when the equation contains an unknown forcing term
while in subsection~\ref{subsec-partial}, we consider the case when
the observations are over a subdomain $[0,L]$ with $L < 1$.

\subsection{Remarks about decay rates in one dimension}
Essentially all the results in the previous section are also
applicable in one dimension with the only difference being the absence
of certain eigenvalues, as we discuss in this subsection. Noting that
$\mu \Delta u + (\lambda + \mu) \nabla (\nabla \cdot u) = \nu u_{xx}$
in one dimension with $\nu = 2\mu + \lambda$, we can simplify
equation~\eqref{erroreq} to the following:
\begin{eqnarray}
   r_t + \rho_0 v_x = - \rho_0 \varphi_\rho \ast v_x \,, \quad
   \rho_0 v_t = \nu v_{xx} -\gamma \rho_0^{\gamma-1} r_x - \rho_0\varphi_u\ast v \,, \label{erroreq1d}
\end{eqnarray}
whereas equation~\eqref{eqwave2} simplifies to
\begin{equation}
  \rho_0 v_{tt}=\nu v_{txx}
+\gamma\rho_0^\gamma v_{xx}+\gamma\rho_0^{\gamma}\varphi_\rho\ast v_{xx}-\rho_0\varphi_u \ast v_t\,.
\label{eqwave2d1}
\end{equation}
From this, we can see that the matrices $M_k$ and $N_k$ from
equations~\eqref{eq:nk} and~\eqref{eq:mk}, respectively, will
remain the same, noting that $k \otimes k = |k|^2 = k^2$. These
$2 \times 2$ matrices will not have eigenvalues $\lambda_{0k}$ or
$\lambda_{dk}$ but only two eigenvalues $\lambda_{\pm k}$ as given in
equation~\eqref{eqeigN} with the same eigenvectors as given before.

Hence the main result in theorem~\ref{mainprop} also applies in
one dimension, with the following modification to the choice of
$\varphi_{uk}$:
\begin{align*}
  \varphi_{uk} = \max\{0;2d-(c_1+c_2)k^2)\} =  \max\left\{0;2d - \frac{\nu 4 \pi^2 k^2}{\rho_0}\right\} \,,
\end{align*}
and the same choice of $\varphi_{\rho k}$ as in
equation~\eqref{varphirho}:
\begin{align*}
 \varphi_{\rho k} = \max\left\{0;\frac{\left((c_1+c_2)|k|^2\right)}{4c_5|k|^2}-1\right\}=\max\left\{0;\frac{\left(\frac{\nu4\pi^2k^2}{\rho_0}+\varphi_{uk}\right)^2}{16\pi^2\gamma\rho_0^{\gamma-1}k^2}-1\right\}.
\end{align*}
Further, in the case of density
observations, the result in theorem~\ref{proprho} is even
stronger, in the sense that with the choice of $\psi_{\rho k}$ and
$\psi_{uk}$ of equations~\eqref{psirho}-\eqref{psiu}, the decay rate
is at least $d$, since there is no non-trivial incompressible flow
field in one dimension.

\subsection{Unknown forcing term}\label{subsec-forcing}

We now consider observers for the linearized one dimensional equations
with a forcing term in the velocity equation. For simplicity, we also
assume that $\rho_0 = 1$, and as earlier, without loss of generality,
we assume $u_0 = 0$.
\begin{eqnarray}
  \rho_t+u_x = 0\,, \quad u_t=\nu u_{xx}-\gamma\rho_x-f (x,t)\,.
\end{eqnarray}
If this forcing term $f(t,x)$ is known, then we also add it to the
velocity equation of the observer~\eqref{eqobslrhou}. Then, by
considering the difference between the observer and original
equations, the forcing term disappears and we are still considering
equation~\eqref{erroreq1d} for the error.

So we now assume that the forcing term is unknown. In this case, we
cannot add it inside the observer equation. So the observer equation
remains unchanged, and then, the difference between the reference
velocity and the observer velocity satisfies
equation~\eqref{eqwave2d1} with a forcing term:
\begin{equation}
  v_{tt}=\nu v_{txx}+\gamma v_{xx}+\gamma\varphi_\rho
  \ast v_{xx} -\varphi_u\ast v_t +f_t \,.
\label{eqwave3}
\end{equation}

We now adapt the spectral analysis. We assume that the mean of $f$ is
$0$ (no bias in the forcing), and then the Fourier decomposition of $f$ is:
$$f(t,x) = \sum_{k\ne 0}f_k(t) e^{i2k\pi x}.$$
Then from~\eqref{eqmode}, we obtain the following new equation for the
$k^{\textrm{th}}$ mode of the error:
\begin{equation}
  a_k''(t) + (\varphi_{uk}+4\nu k^2\pi^2)\, a_k'(t) + 4k^2\pi^2\gamma
  (1 + \varphi_{\rho k}) \, a_k(t) = f_k'(t).
\label{eqmode2}
\end{equation}
We just need to find a particular solution to this equation, and add
it to the general solution that we found in section~\ref{subsec-spectral}.

We consider a very simple case (for clarity reasons and for a better understanding of the phenomenon), where the time dependence of the forcing is
a sine (or cosine) function:
\begin{equation}
  f_k(t) = c_k \sin(2\omega_k \pi t),
\end{equation}
where $\omega_k$ is the frequency of the forcing oscillation of mode
$k$. Defining $\alpha_k=\varphi_{uk}+4\nu k^2\pi^2$, and
$\beta_k=4k^2\pi^2\gamma(1+\varphi_{\rho k})$,
equation~\eqref{eqmode2} becomes:
\begin{equation}
  a_k''(t)+\alpha_k a_k'(t) + \beta_k a_k(t) = 2c_k\omega_k\pi
  \cos(2\omega_k\pi t).
\end{equation}
A particular solution is then given by
\begin{equation}
a_k(t) = A_k \cos(2\omega_k\pi t) + B_k \sin(2\omega_k\pi t).
\end{equation}
The constants $A_k$ and $B_k$ are solution of the following linear system:
$$ \begin{array}{l}
  A_k(\beta_k-4\omega_k^2\pi^2)+B_k(2\omega_k\pi\alpha_k) =
  2c_k\omega_k\pi,\\ A_k(-2\omega_k\pi\alpha_k) +
  B_k(\beta_k-4\omega_k^2\pi^2) = 0.
\end{array} $$
Then, we get:
\begin{eqnarray}
  A_k &=& 2c_k \omega_k \pi \frac{\beta_k-4\omega_k^2\pi^2}
  {(\beta_k-4\omega_k^2\pi^2)^2 + (2\omega_k\pi\alpha_k)^2} , \\ B_k
  &=& 2c_k \omega_k \pi \frac{2\omega_k\pi\alpha_k}
  {(\beta_k-4\omega_k^2\pi^2)^2 + (2\omega_k\pi\alpha_k)^2}.
\end{eqnarray}
The amplitude of the particular solution is then given by
\begin{equation}
  D_k = \sqrt{A_k^2+B_k^2} = \frac{2c_k\omega_k\pi}
  {\sqrt{(\beta_k-4\omega_k^2\pi^2)^2 + (2\omega_k\pi\alpha_k)^2}}
\end{equation}

Increasing $\varphi_{\rho k}$ or $\varphi_{u k}$ (or both) will
make $\beta_k$ or $\alpha_k$ (or both) increase, and then $D_k$ will
decrease. This means that we can make $D_k$ become as small as we
want and the observer will converge towards the true state. 
But of course, the numerical performance of the observer is severely
degraded in comparison with previous cases, as it is usually not
possible to consider extremely high values of feedback coefficients
from a numerical point of view.

Concerning the density, adapting the previous calculations on the
decrease of $\rho$ (knowing the decrease of $u$), we get the following amplitude of the particular solution:
\begin{equation}
  E_k = \frac{k}{\omega_k} (1+\varphi_{\rho k})\, D_k,
\end{equation}
which means that one can adapt the limit amplitude of $\rho$ by
changing the values of $\varphi_{\rho k}$. Theoretically, choosing
$\varphi_{\rho k} = -1$ has the effect of completely removing the influence of the forcing term on
the density, but of course, it is not numerically stable or
physically consistent to consider negative feedback coefficients.

\subsection{Heuristic in the case of observations over a part of the domain} \label{subsec-partial}

In the case of observations over only a part of the domain, we define
the observation domain $\Omega = [0,L]$ with $L < 1$. In this case, in
order to proceed with Fourier analysis, we will need to find the
Fourier transform of $v \mathbbm{1}_\Omega$ because the feedback terms
$\varphi\ast v$ will be replaced with
$\varphi\ast (v \mathbbm{1}_\Omega)$. The Fourier series for
$\mathbbm{1}_\Omega$ is
\begin{equation}\label{fourier1}
  \mathbbm{1}_\Omega = L + \sum_{k=1}^\infty \left[ \frac{\sin 2\pi
      kL}{k} \sin 2\pi k x + \frac{1-\cos 2\pi kL}{k} \cos 2\pi k x
    \right] \,.
\end{equation}
Thus we see that the Fourier series of $v \mathbbm{1}_\Omega$
will have Fourier components for all $k$ even in the case when $v$ has just a single Fourier mode. 

Such coupling of Fourier modes will numerically slightly degrade the
performance of the observer, as at any time, some energy will be
transferred between different Fourier modes, as in the full nonlinear
model.

If we rewrite equation~\eqref{eqwave2d1} (for simplicity with
$\rho_0=1$), we get:
\begin{equation}
  v_{tt}=\nu v_{txx}
+\gamma v_{xx}+\gamma\varphi_\rho\ast \left(v_{xx}\mathbbm{1}_\Omega\right) - \varphi_u \ast \left(v_t\mathbbm{1}_\Omega\right)\,.
\end{equation}
Then, for simplicity reasons, assuming $v$ only has one single
Fourier mode $k$, equation~\eqref{eqmode2} rewrites:
\begin{equation}
  a_k''(t) + (L \varphi_{uk}+4\nu k^2\pi^2)\, a_k'(t) +
  4k^2\pi^2\gamma(1+L \varphi_{\rho k})\,a_k(t)= 0.
\end{equation}
Indeed, only the 0th order term ($L$, from~\eqref{fourier1}) in the Fourier decomposition of
$\mathbbm{1}_\Omega$ will be kept through the convolution with 
$v$ (or one of its derivatives). So, the decay rate becomes
\begin{equation}
\frac{\varphi_{uk}L+4\nu k^2\pi^2}{2}.
\label{eqeigL}
\end{equation}
Of course, this is an approximation, as even if $v$ only has a
single Fourier mode at time $t=0$, the convolution with a
characteristic function leads to mode mixing, as in a nonlinear
situation. But we assume here that most of the energy is along
mode $k$ (if only this mode is present at the initial time), which
will be confirmed by numerical experiments in next section.

\section{Numerical experiments on the 1D compressible Navier-Stokes
  observer} \label{sec-numeric}

In this section, we report some of the numerical investigations in one dimension that
illustrate the linear theory we discussed above. We also present
results of using the same observer as in the linear case for two
additional scenarios,
namely, (i) the observer when the velocity is observed only over a
subinterval of the domain in section~\ref{subsec-partial-num}, and
(ii) in section~\ref{subsec-nonlin-num}, observer for the fully
nonlinear equations.

\subsection{Numerical configuration} \label{subsec-code}

The space domain is $[0,1]$ and we still assume periodic boundary conditions. The
discretization involves $10^2$ grid points, with a step $\Delta
x=10^{-2}$. We consider a time step $\Delta t=10^{-3}$. We also
experimented with increasing the spatial resolution (and
correspondingly decreasing the time step), but the results are almost
identical and not presented here. The adiabatic exponent is set to
$\gamma = 1.4$, and the diffusion is set to $\nu=5.10^{-2}$ (except in
section~\ref{subsec-nonlin-num}). The numerical code uses a conservation
form of the compressible Navier-Stokes system, with $\rho$ and
$\rho u$ as variables. A finite volume scheme is used, in which the
inviscid flux is computed using an approximate Riemann solver (e.g.
VFRoe scheme). Time integration scheme is a third order explicit
Runge-Kutta scheme, where the time step is chosen based on a CFL
condition. 

In order to reproduce a quasi-linear situation, we consider the true
(observed) solution $\rho(t,x)= 1$ and $u(t,x) = 0$ while the initial
conditions for the observer are set to:
\begin{align}
\begin{array}{l} \hat\rho_I(x)= 1+5.10^{-2}\,\sin(2\pi k x), \ \ 
  \hat u_I(x) = 5.10^{-2}\,\sin(2\pi k x),
\end{array} \label{eqlinpert}  
\end{align}
so that the mean values of $\hat\rho$ and $\hat u$ are $\rho_I=1$ and
$u_I=0$ respectively, and where $k$ is a given mode, usually the first
one ($k=1$, unless differently specified).
\begin{figure}[t!]
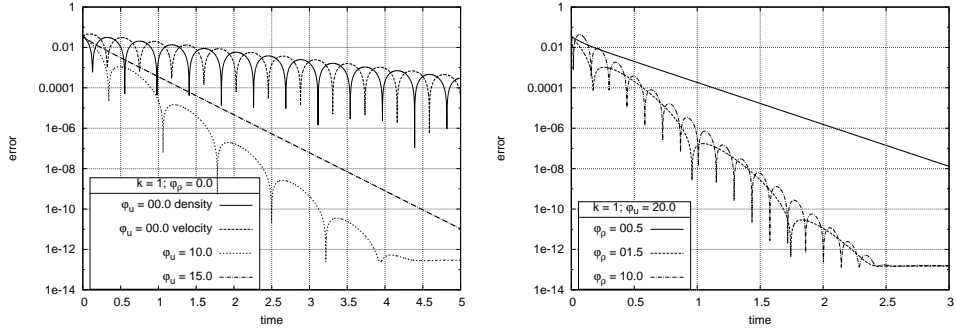

  \centering
  \includegraphics[width=0.49\textwidth]{{{figures/obs001_pr00.0_new}}}
  \includegraphics[width=0.49\textwidth]{{{figures/obs001_pu20.0_new}}}
  \caption{The $L^2$ norm of the difference between the observer
    $(\hat\rho,\hat u)$ and the solution $(\rho,u)$ versus time. Solid
    and dotted lines are the errors in $\rho$ and $u$,
    respectively. The left panel is for fixed $\varphi_\rho = 0$ with
    varying $\varphi_u$ while the right panel is for fixed $\varphi_u = 20$
    with varying $\varphi_\rho$.}
  \label{fig00}
\end{figure}

We first look at the solution without any feedback
term, see equation~\eqref{eqobslrhou} with $\varphi_u=\varphi_\rho=0$. Figure~\ref{fig00} (solid and dashed curves, left panel) shows
the evolution (in log scale) of the $L^2$ norm of the difference
between the observer $(\hat\rho,\hat u)$ and the solution
$(\rho,u)$. As there is no feedback, all the Fourier coefficients
$\varphi_{\rho k}$ and $\varphi_{uk}$ are equal to $0$, and then from
equation~\eqref{eqdelta}, the discriminant is $\Delta\simeq -217.2$,
and the theoretical decay rate (only due to diffusion) is given by
equation~\eqref{eqeigN}: $d_{th}=0.987$.  Also the oscillation period
can be computed from~\eqref{eqeigN} and~\eqref{eqdelta}:
$\omega_{th}=\displaystyle\frac{4\pi}{\sqrt{-\Delta}}\simeq 0.85$.

Numerically, the slope of the solid curve in the left panel of
figure~\ref{fig00} gives a numerical decay rate
$d_{\textrm{num}}=0.980$. Note that the figures show the errors to
base 10, hence the slope of the semilog plot is $d_\textrm{{num}} /
\log(10) = 0.426$.  The numerical oscillation period is approximately
$\omega_\textrm{{num}}=0.852$. Note that one period corresponds to two
oscillations on the figure for the norm of the cosine. This excellent
agreement between theoretical and numerical values can be reproduced
for other modes and other values of the parameters.

\subsection{Simple nudging observer} \label{subsec-nudge}

We now consider the nudging framework (see section~\ref{subsec-remark} and equation~\eqref{erroreqnudg1}), and we hence suppose that there is some feedback only in the velocity
equation. We then first let $\varphi_\rho=0$, and only modify the
values of $\varphi_u$. This simulates the nudging, or asymptotic
observer: as only the velocity is measured, only the velocity is
corrected in the observer system. Table~\ref{tab1} shows the
theoretical and numerical decay rates and oscillation periods for
several values of $\varphi_{u1}$.
\begin{table}[t!]
\centering
$
\begin{array}{c|cc|cc}
\varphi_{u1} & \textrm{\begin{tabular}{c}Theoretical \\ decay rate\end{tabular}}
& \textrm{\begin{tabular}{c}Numerical\\ decay rate \end{tabular}}
& \textrm{\begin{tabular}{c}Theoretical\\oscillation period \end{tabular}}
& \textrm{\begin{tabular}{c}Numerical\\ oscillation period\end{tabular}} \\
\hline
0   & 0.987 & 0.980 (0.426) & 0.85 & 0.86\\
0.1 & 1.037 & 1.032 & 0.85 & 0.85\\
0.5 & 1.237 & 1.237 & 0.86 & 0.85\\
1   & 1.487 & 1.486 & 0.86 & 0.87\\
5   & 3.487 & 3.485 & 0.97 & 0.98\\
10  & 5.987 & 6.012 (2.61) & 1.42 & 1.38\\
12.895 & 7.434 & 6.590 & 81.0 & - \\
15  & 4.393 & 4.364 (1.895) & - & - \\
20  & 2.897 & 2.861 & - & - \\
\end{array}
$
\caption{Theoretical and numerical decay rates and oscillation periods
  for several values of $\varphi_u$ ($k=1$, $\varphi_\rho=0$). The
  values in parenthesis give the numerical decay rates in base-10, in
  order to compare with slopes of lines in figure~\ref{fig00}}
\label{tab1}
\end{table}

The first remark is that the numerical results perfectly match the
theoretical results, except for the particular value of
$\varphi_{u1}=12.895$. In this case, the numerical decay rate is
slightly smaller than the theoretical one. Also, no oscillations can
be seen on the results, which is reasonably in agreement with a
theoretical period of $81$ which will be impossible to see with a
final time of $T = 5$.

As $\varphi_{u1}$ increases, the decay rate increases, until
$\varphi_{u1}$ reaches $4\pi\sqrt{\gamma}-4\nu\pi^2 \simeq 12.895$
(see remark in section~\ref{subsec-remark}), for which the discriminant is
equal to $0$, and then positive for increasing $\varphi_{u1}$. The
corresponding optimal decay rate is $d_{th}=2\pi\sqrt{\gamma}\simeq
7.434$ (see section~\ref{subsec-remark}). We can see that the decay rate
then decreases, as the discriminant takes larger positive values, so
that one of the two eigenvalues gets closer to $0$.

One can see on figure~\ref{fig00} (left panel, dotted and dash-dotted
curves) that the error decreases much stronger than the case of no
feedback (figure~\ref{fig00}, left panel, solid curve). We can also see
that with increasing $\varphi_u$, the period of oscillations increases
and eventually there are no oscillations (discriminant is
positive). It confirms that the decay rate decreases if $\varphi_{u1}$
is increased too much.

Similar results have been observed for other values of the diffusion
$\nu$, and for other modes $k$, and in each case, these numerical
results match well with the theoretical predictions.

\subsection{Results on the full state observer} \label{subsec-full}

We now use the full state observer, with an additional feedback term in the
density equation (see equation~\eqref{erroreq2}). We refer here to theoretical results from sections~\ref{subsec-spectral} and~\ref{subsec-result}. We first set $\varphi_{u1}=20$, for which the
discriminant is positive, the largest eigenvalue gets closer to $0$,
and then the decay rate becomes non optimal.

\begin{table}
\centering
$
\begin{array}{c|cc|cc}
\varphi_{\rho 1} & \textrm{\begin{tabular}{c}Theoretical \\ decay rate\end{tabular}}
& \textrm{\begin{tabular}{c}Numerical\\ decay rate \end{tabular}}
& \textrm{\begin{tabular}{c}Theoretical\\oscillation period \end{tabular}}
& \textrm{\begin{tabular}{c}Numerical\\ oscillation period\end{tabular}} \\
\hline
0   & 2.897 & 2.861 & - & - \\
0.5   & 4.838 & 4.790 (2.080) & - & - \\
1.184 & 10.94 & 9.74 &  &  \\
1.5 & 10.98 & 11.03 (4.791) & 1.50 & 1.57 \\
5   & 10.99 & 11.03 & 0.43 & 0.43 \\
10  & 10.99 & 11.06 (4.802) & 0.28 & 0.28 \\
\end{array}
$
\caption{Theoretical and numerical decay rates and oscillation periods
  for several values of $\varphi_{\rho}$ with fixed $\varphi_u=20$ for
  the $k=1$ mode. (The values in parenthesis are again decay rates in
  base-10 for comparison with figure~\ref{fig00}.)}
\label{tab2}
\end{table}
Table~\ref{tab2} shows the theoretical and numerical decay rates (and
oscillation periods) for several values of $\varphi_{\rho 1}$.  As the
discriminant is positive when $\varphi_{u1}=20$ and $\varphi_{\rho
  1}=0$, equation~\eqref{varphirho} gives the theoretical value for
which the discriminant comes back to negative values: $\varphi_{\rho
  1}=\displaystyle\frac{(\varphi_{u1}+4\nu\pi^2)^2}{16\gamma\pi^2}-1\simeq 1.184$.

We observe numerically that the decay rate is very similar to the
theoretical rate, and we clearly observe the transition from positive
to negative discriminant with the stabilization of the decay rate, and
the apparition of oscillations. Increasing $\varphi_{\rho 1}$ to a
much larger value than the optimum given by~\eqref{varphirho} is not
necessary, as the decay rate does not increase, and the period of
oscillations increases quite quickly. This is clearly seen from the
plots in the right panel of figure~\ref{fig00}.

As we have seen, by adding the derivative of the velocity as a
feedback to the density equation, we were able to significantly
increase the decay rate of the error, in comparison with only a
feedback in the velocity equation.

\subsection{Observers with observations over a part of the domain} \label{subsec-partial-num}

\begin{figure}[t!] \begin{center}
    \includegraphics[width=0.49\textwidth]{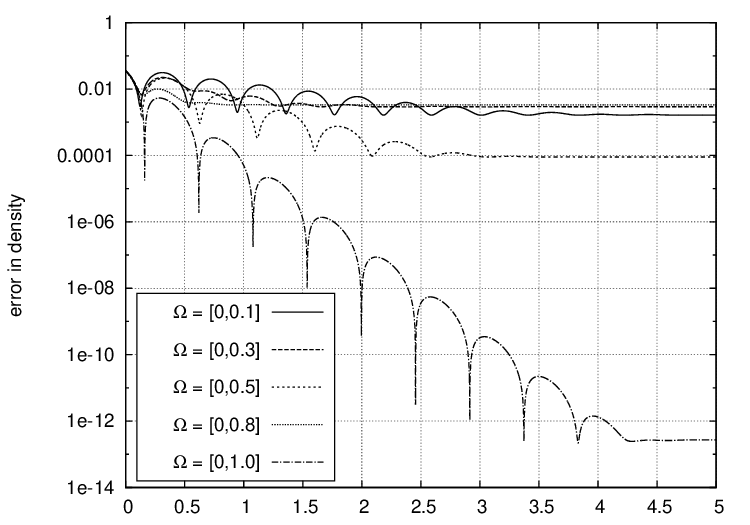}
    \includegraphics[width=0.49\textwidth]{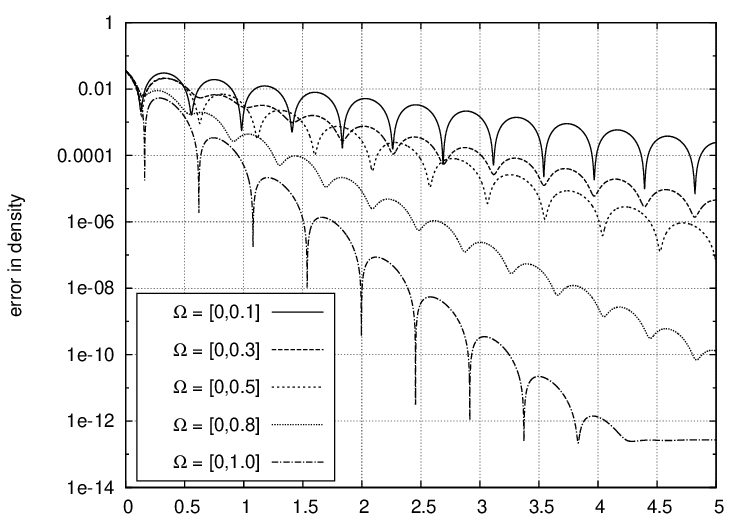}
    \includegraphics[width=0.49\textwidth]{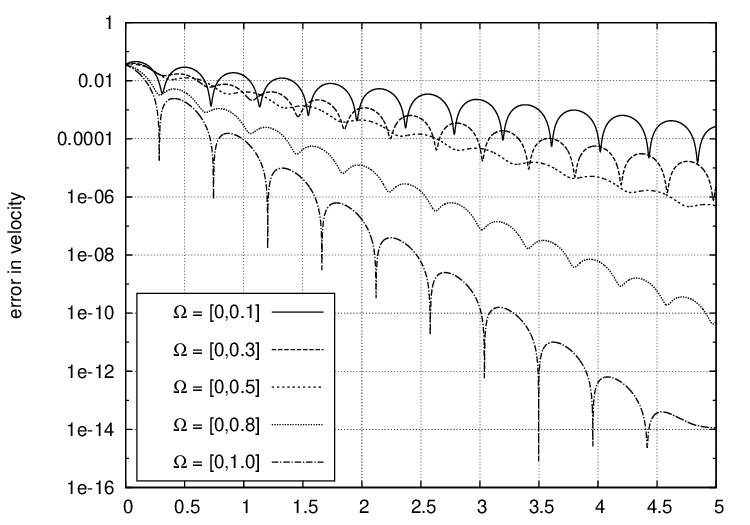}
    \includegraphics[width=0.49\textwidth]{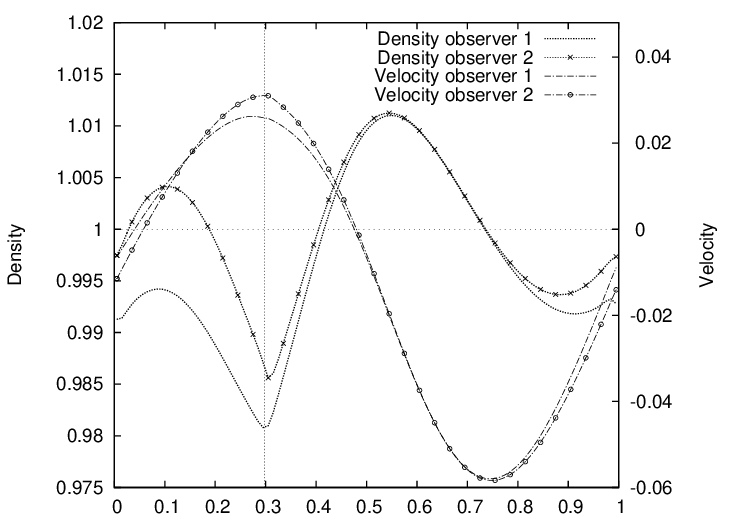}
    \caption{$L^2$ norm of the difference between the observer
      $\hat\rho, \hat u$ and the solution $\rho,u$, in the case of
      observations over the sub-interval $[0,L]$ for $L = 0.1, 0.3,
      0.5, 0.8, 1$, all of them with strength of feedback terms being
      $\varphi_\rho = 0.5$ and $\varphi_u = 10$. Left panels are for
      feedback from equation~\eqref{feedbackrhou} while the top right
      panel is for feedback from equation~\eqref{newfeedbackrho}. The
      error in velocity for the feedback from
      equation~\eqref{newfeedbackrho} is almost identical to the
      bottom left panel and hence not shown. The bottom right panel
      shows the observer solutions at time $t = 0.16$. Note that
      ``observer 1'' refers to~\eqref{feedbackrhou} while ``observer
      2'' refers to~\eqref{newfeedbackrho}.}
    \label{partfig}
\end{center} \end{figure}

In this section we present the numerical results of the full state observer,
see equation~\eqref{feedbackrhou}, but with observations over only a part of
the domain. The left panels
of figure~\ref{partfig} shows the decay rate of the $L^2$ error
between the observer $(\hat{\rho}, \hat{u})$ and the actual solution
$(\rho,u)$.

We see that, as expected from section~\ref{subsec-partial}, the rate of decay is smaller for smaller
observational intervals. We also see that the error in velocity
decreases linearly (with oscillations) whereas the error in density
saturates at a fairly high but constant value. This is because the
solution of the observer equation converges to a solution with
$\hat{u} = 0$ but with $\hat{\rho} = \rho_n$ where $\rho_n \ne \rho_I$
- i.e. the observer density is shifted by an amount which compensates
for the initial discrepancy between the mean of $\hat{\rho}$ and the
mean of $\rho$. This is not surprising since the original equations
themselves are invariant under the constant shift in density.

In order to overcome this problem in the case of observations over partial domain,
we propose the following modification of the feedback terms in
equation~\eqref{feedbackrhou}:
\begin{eqnarray}
  F_\rho(\hat\rho,\hat u,u) &=& \varphi_\rho\ast D_\rho \left[
    \mathbbm{1}_\Omega (u-\hat u) - \langle u-\hat u \rangle
    \right], \label{newfeedbackrho}
\end{eqnarray}
where $\langle f \rangle$ indicates average of $f$ over the interval
$[0,L]$. This ensures that the average of the feedback term is zero
and hence the equilibrium solution of this equation also has mean
zero. Note that in the case of $L = 1$, i.e., the case of full
observations, this average is just zero and the feedback in
equation~\eqref{newfeedbackrho} is identical to that in
equation~\eqref{feedbackrhou}. The errors obtained by using this new
observer are shown in the right upper panel of
figure~\ref{partfig}. We clearly see that the observer $\hat\rho$ now
approaches the solution $\rho$ and the error decreases as expected.

In order to clearly see the effect of the observer, the lower right
panel of figure~\ref{partfig} shows the actual observer solutions for
the case when $\Omega = [0, 0.3]$. They clearly show the effect of
incorporating the observations, and also the difference between the
observer with the incorrect mean and the one with correct mean. The
effect is of course more pronounced on the density than on the
velocity.

\begin{figure}[t!] \begin{center}
    \includegraphics[width=0.49\textwidth]{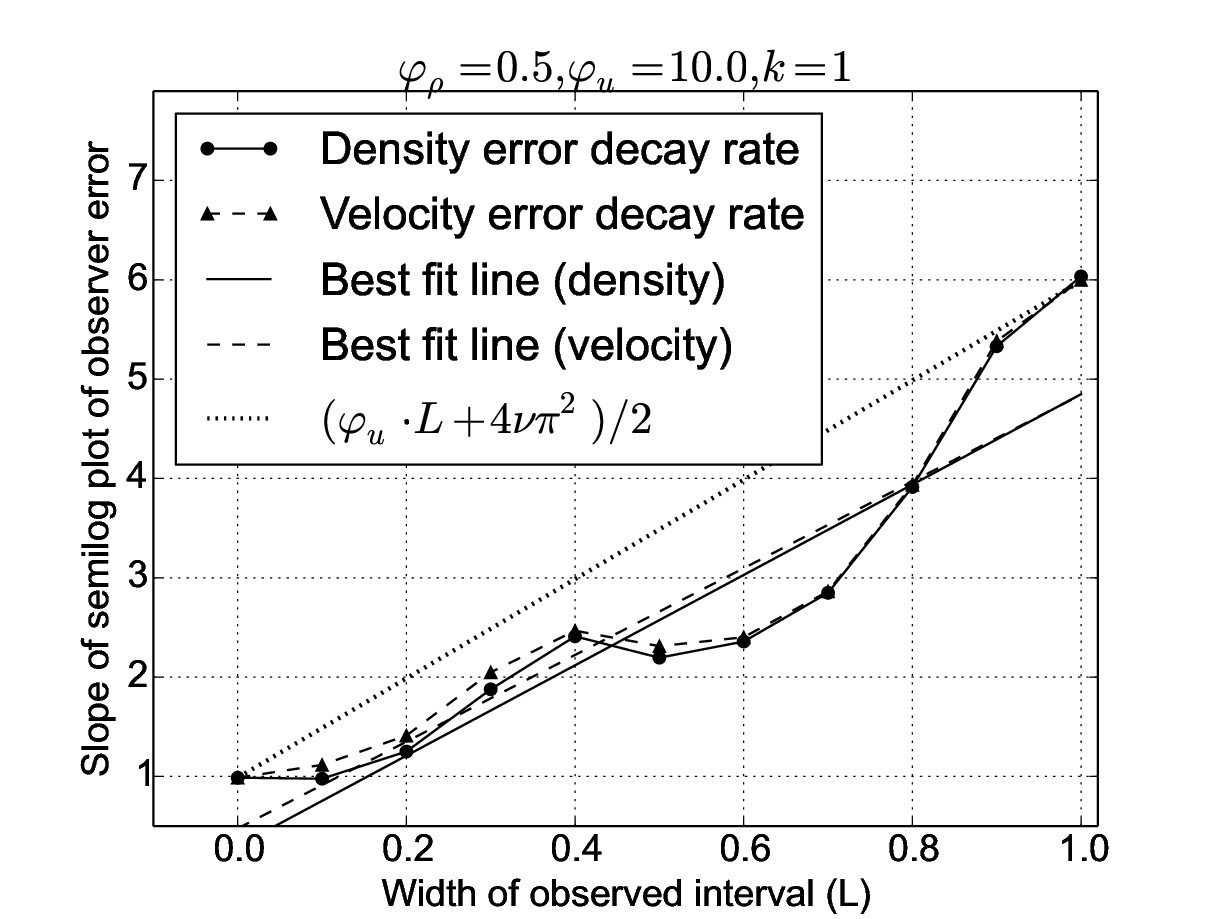}
    \includegraphics[width=0.49\textwidth]{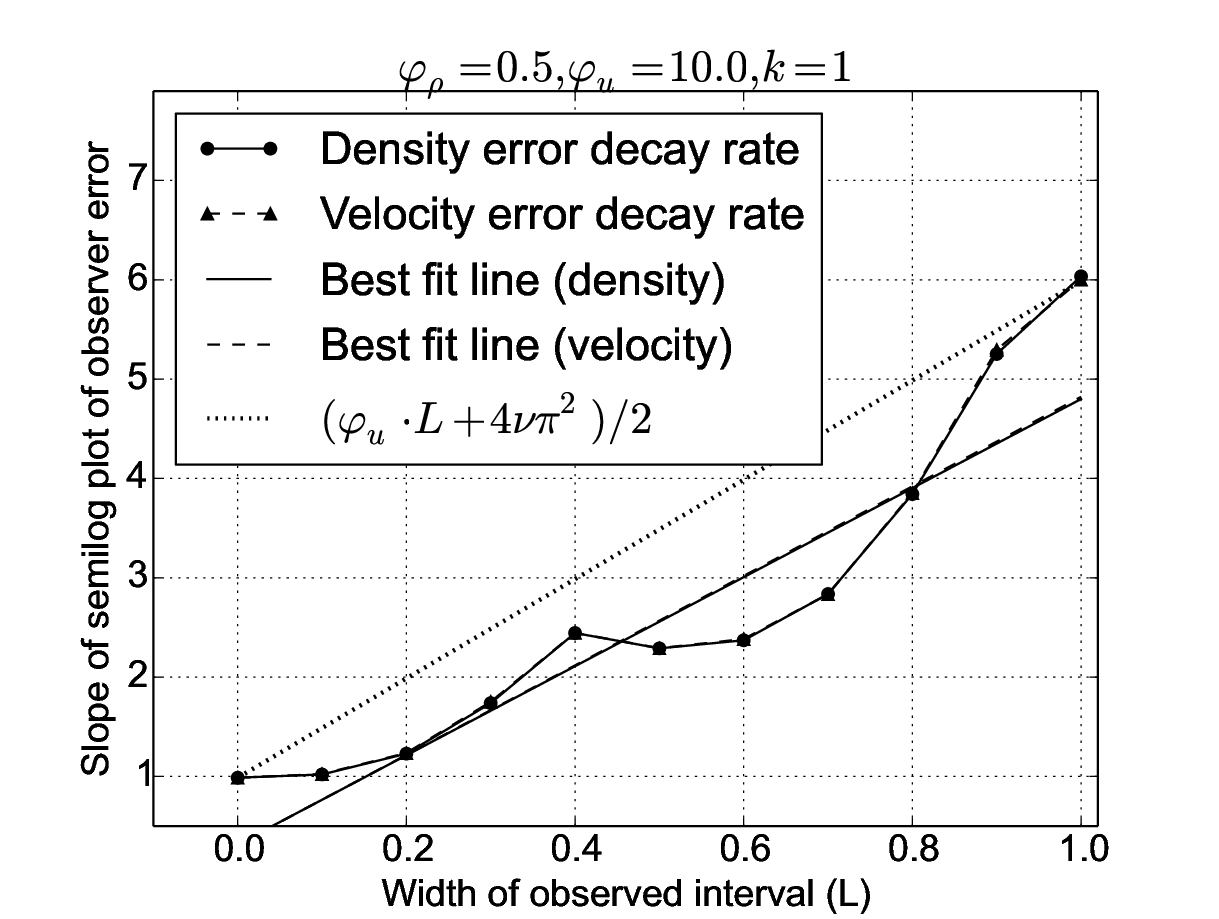}
    \caption{The decay rates for the velocity and the density, as a
      function of the length of the observation interval (quasi-linear
      case in the left panel and fully nonlinear case in the right
      panel).  Note that the slope of the best fit line $\approx 4.5$
      whereas the slope of best fit line by taking only the points for
      $L = 0.1, 0.2, 0.3, 0.4, 0.9, 1.0$ is $\approx 5.7$.}
    \label{partobsrate}
\end{center} \end{figure}

In this case the discriminant of the decay rate is clearly negative
(as evidenced by oscillations in the plot for errors) but we can
calculate the decay rate.  Figure~\ref{partobsrate} shows the decay
rate of the observer as a function of the length of the interval over
which the velocity is observed, for the case when $\varphi_\rho = 0.5$,
$\varphi_u =10$ for the $k=1$ mode. We notice that this is pretty close
to a straight line and the best fit line is given as follows:
\begin{align}
  D = 4.99 L  + 0.24 \,.
\end{align}
Comparing this with the rate given in equation~\eqref{eqeigL}, we see
that the last equality is a reasonable approximation. Thus even though
we cannot calculate the exact decay rate (see 
section~\ref{subsec-partial}), a reasonable estimate can be
obtained by using formula~\eqref{eqeigL}:
\begin{align}
  D \approx \frac{\varphi_u L + 4 \nu \pi^2 k^2}{2} = 5 L + 0.987
\end{align}
in experiments corresponding to figure~\ref{partobsrate}
($\varphi_u=10$).

This figure shows that the global slope of the numerical decay rates
versus the length of the observation interval is close to 
$\frac{\varphi_u}{2}=5$. There is also a good agreement between
numerical decay rates and the approximated theoretical ones (see
equation~\eqref{eqeigL}) when $L$ is either small (almost no observations)
or large (almost all the domain is observed), while the numerical
decay rates are degraded when only half the domain is observed.
When $L\approx 0.5$, previous spectral studies show that mixing of 
Fourier modes has a higher effect than when $L\approx 0$ or $1$.
Convolution with the characteristic function leads to non negligible
transfers of energy from mode $k$ to other modes, and then to a 
smaller decay rate.

But figure~\ref{partobsrate} also shows that even with a small
observed subdomain, our observer is still very efficient and we can
still control the decay rate by increasing $\varphi_u$. We note that
we have not considered the problem of placement of the observations,
which is certainly important and has a huge influence on the decay
rate \cite{privat2015optimal, hebrard2003optimal}.

\subsection{Observer for nonlinear Navier-Stokes system} \label{subsec-nonlin-num}

In this section, we will report the numerical results of using the
full state observer for the nonlinear system of equations, i.e. from
equations~\eqref{eqobsrhou}, with the feedback terms as in
equations~\eqref{feedbackrhou}, in one dimensional system. Note that we do not have theoretical
estimates of the decay rates towards the equilibrium solution. The aim
of this section is to understand the efficiency of the above nonlinear
observer.

\begin{figure}[t!] \begin{center}
    \includegraphics[width=0.49\textwidth]{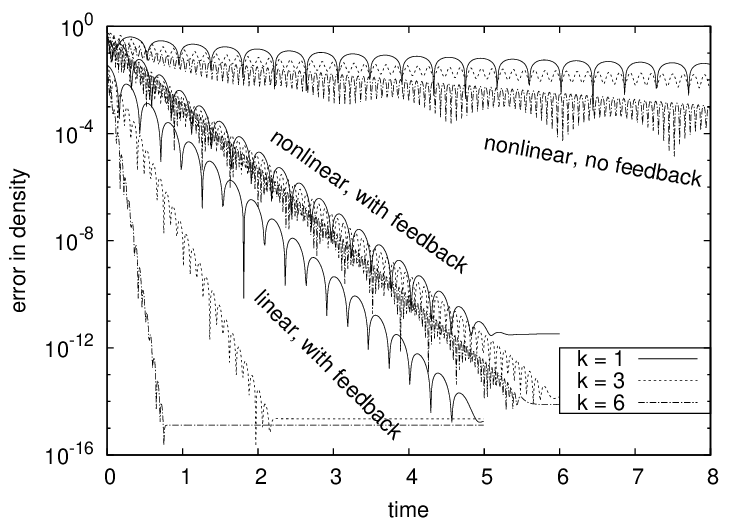}
    \includegraphics[width=0.49\textwidth]{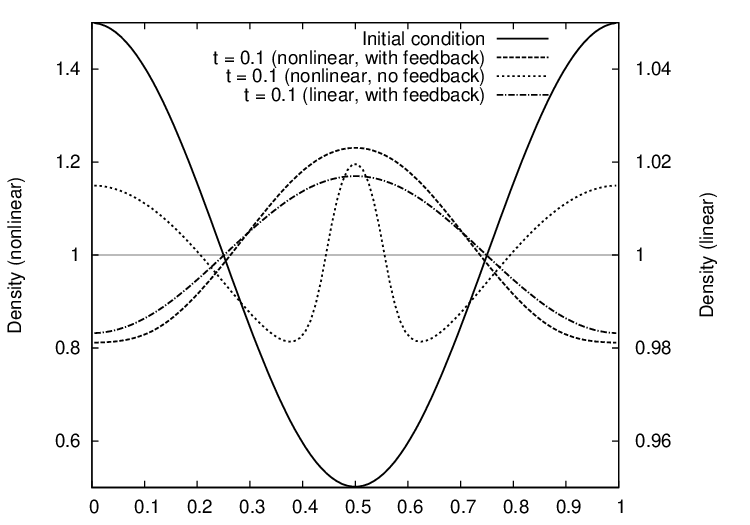}
    \includegraphics[width=0.49\textwidth]{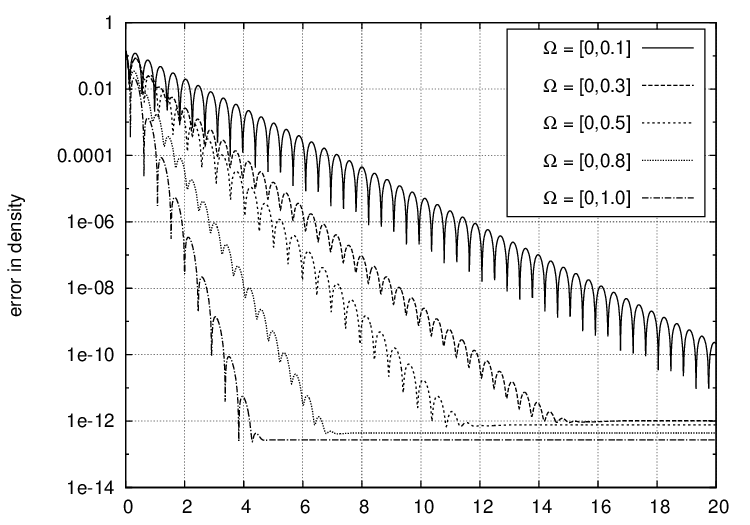}
    \includegraphics[width=0.49\textwidth]{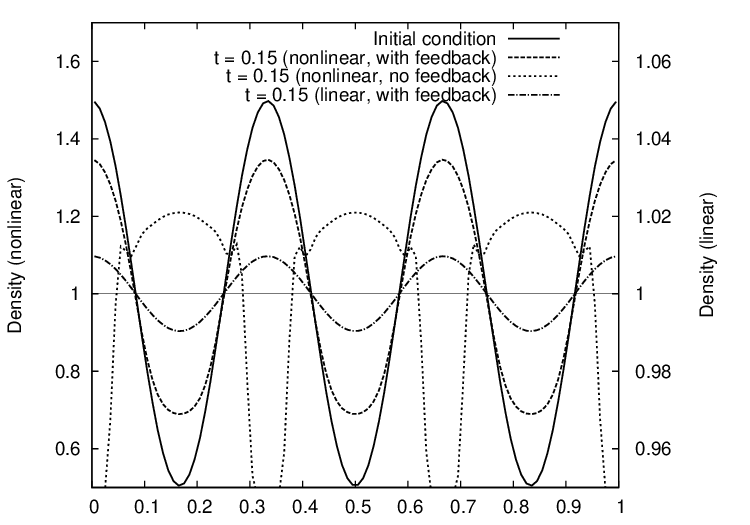}
    \caption{The decay rates of the nonlinear observer compared with
      that of the linear observer (top left). The actual observer
      solutions are also shown ($k=1$ in the top right, $k=3$ in the
      bottom right). The decay rates for the nonlinear observer with
      observations over partial domain is in the bottom left.}
    \label{fignonlin}
\end{center} \end{figure}

Firstly, we consider observations of the equilibrium solution, so that
we are essentially studying the decay of the observer towards the
equilibrium. The left panel of figure~\ref{fignonlin} shows the decay
rates of the observer solution with the following two initial
conditions, with the feedback term set to zero, and for the case with
$(\varphi_\rho, \varphi_u) = (0.2, 10)$:
$$ \begin{array}{l} \hat\rho_I(x)= 1+5.10^{-1}\,\sin(2\pi k x),\quad \hat
  u_I(x) = 5.10^{-1}\,\sin(2\pi k x),
\end{array} $$
Note that the perturbation strength is 10 times larger than the
``quasi-linear'' case with the perturbation in
equation~\eqref{eqlinpert}. We consider two cases, $k=1$ and
$k=3$. The cases with other values of $k > 1$ show a very
similar behavior to the case $k=3$ and hence is not shown here explicitly. For
reference, we have also plotted the decay rates of the linear case as
well.

The right panel of that figure shows the actual observer solutions for
the $k = 1$ case. Note the different scales for the nonlinear (left
axis) and linear (right axis) regimes. We clearly see that with
perturbation strength of $5.10^{-1}$, even though the initial
condition only has $k=1$ mode present, at time $t=0.1$, higher modes
are excited (dotted line, with no feedback). When the observer
feedback is added, the solution very quickly decays before higher
modes are excited, bringing the perturbation to the level where the
linear theory is a good approximation and thus the decay rate is
identical to that of the linear case of perturbation strength
$5.10^{-2}$.

The case for $k = 3$ (and indeed all higher modes with $k > 1$) is
quite interesting. In this case, the initial condition of $k = 3$
excites the $k = 1$ mode. Thus even when the observer feedback is
added, the decay rate is not the same as the linear $k=3$ decay rate
but rather it is closer to the linear $k = 1$ rate. The same behavior
is seen for other modes with $k > 1$. Figure~\ref{fignonlin} shows
example of this decay for $k=3$ and $k=6$, along with the actual observers. It
is quite clear that in the quasi-linear case, modes other than the one
contained in the initial condition are not excited while in the
nonlinear case, they are quite clearly excited.

Finally, we also performed numerical experiments, still with the full state observer on the nonlinear system, but now with partial observations over various domain sizes, as
discussed in section~\ref{subsec-partial}. Quite surprisingly, the
behavior in the nonlinear case is very similar to the linear case: the
decay rate is close to being linear in the size of the domain - see
the right panel of figure~\ref{partobsrate} and the bottom left panel
of figure~\ref{fignonlin}.

\section{Conclusion} \label{sec-conclude}

In this paper, we were interested in observer design for a viscous
compressible Navier-Stokes equations. Thanks to intrinsic properties of the system
(symmetries), we were able to design observers in order to reconstruct
the full solution (both velocity and density) when only one variable
(either velocity or density) is observed, even partially.

A spectral study showed that for a tangent linear system (linearized
around the equilibrium state of constant velocity and density), we can prove the convergence of the
observer towards the solution, and we can control the decay rate
of the error, with explicit formulas for the feedback coefficients
as functions of the desired decay rate.

Numerical experiments in one dimension are in perfect agreement with theory in the
linearized situation: we can obtain any decay rate by increasing the
observer coefficients. Numerical experiments also show that our
observer is still very efficient in the case of observations over partial domain,
and also in the full nonlinear case.

The application of this kind of observer in other specific cases  such as
full primitive equations of the
ocean or the atmosphere, will be a natural extension
of this work.  Observers for fluid equations coupled to
other quantities such as temperature or salinity with observations of
these fields instead of velocity observations will also be an
interesting extension that will be of great interest in practice since
these types of measurements are more common than measurements of the velocity field. It will also be quite
challenging and interesting for practical applications to consider
observers in the case when observations are available discretely in
time and/or in space.

\bibliographystyle{siam}
\bibliography{observer-cns}

\end{document}